\documentclass[12pt,twoside]{amsart}
\usepackage{amssymb}
\usepackage{amscd}
\usepackage{color}

\usepackage{tikz}
\usetikzlibrary{cd}

\usepackage{hyperref}

\title[Relative semi-ampleness in positive characteristic]
{Relative semi-ampleness in positive characteristic}  

\author{Paolo Cascini and Hiromu Tanaka}
\address{Department of Mathematics, Imperial College London, 180 Queen's Gate, 
London SW7 2AZ, UK} 
\email{p.cascini@imperial.ac.uk}
\address{Graduate School of Mathematical Sciences, 
The University of Tokyo, 
3-8-1 Komaba, Meguro-ku, Tokyo 153-8914, JAPAN} 
\email{tanaka@ms.u-tokyo.ac.jp}
\thanks{The first author was funded by EPSRC. 
The second author was funded by EPSRC and 
the Grant-in-Aid for Scientific Research (KAKENHI No. 18K13386). 
We would like to thank Y. Gongyo, Z. Patakfalvi and S. Takagi for many useful discussions. We would also like to thank B. Bhatt for informing us about \cite{BS17}. 
We would like to thank the referee for many useful comments.}

\newcommand{\alt}[0]{{\operatorname{alt}}}

\newcommand{\codim}[0]{{\operatorname{codim}}}
\newcommand{\red}[0]{{\operatorname{red}}}
\newcommand{\Ker}[0]{{\operatorname{Ker}}}

\newcommand{\Proj}[0]{{\operatorname{Proj}}}
\newcommand{\Spec}[0]{{\operatorname{Spec}}}

\newcommand{\Supp}[0]{{\operatorname{Supp}}}
\newcommand{\Pic}[0]{{\operatorname{Pic}}}

\newtheorem{thm}{Theorem}[section]
\newtheorem{lemma}[thm]{Lemma}

\newtheorem{prop}[thm]{Proposition}

\theoremstyle{definition}
\newtheorem{ex}[thm]{Example}

\newtheorem{dfn}[thm]{Definition}

\newtheorem{rem}[thm]{Remark}

\newtheorem{step}{Step}
\newtheorem*{claim}{Claim}

\makeatletter
  
  \@addtoreset{equation}{section}
  \makeatother

\newcommand{\MO}{\mathcal{O}}
\newcommand{\F}{\mathbb{F}}

\newcommand{\Q}{\mathbb{Q}}
\newcommand{\Z}{\mathbb{Z}}
\newcommand{\p}{\mathfrak{p}}
\newcommand{\q}{\mathfrak{q}}
\newcommand{\m}{\mathfrak{m}}

\newtheorem{theorema}{Theorem}

\begin{document}

\maketitle

\begin{abstract}
Given an invertible sheaf on a fibre space between projective varieties of positive characteristic, we show that  fibrewise semi-ampleness implies  relative semi-ampleness. 
The same statement fails in characteristic zero. 
\end{abstract}

\tableofcontents

\setcounter{section}{0}

\section{Introduction}
It is a fundamental 
problem in algebraic geometry to study under what conditions a nef line bundle on a projective variety  is semi-ample. 
For instance, the abundance conjecture predicts that, on a minimal projective variety,  
the canonical divisor is always semi-ample. 
On the other hand,  it is not easy
in general to find criteria that hold for any line bundle. 

Over a field of positive characteristic, 
it seems that semi-ampleness sometimes behaves better 
than in characteristic zero.  
One of the most typical examples is Keel's result \cite{keel99}, 
which has recently played a crucial 
role in the minimal model program of positive characteristic (e.g. see \cite{hx13}). 

 The goal of this paper is to provide 
a necessary and sufficient condition under which, given a  morphism of $\mathbb F_p$-schemes $f\colon X\to Y$,  an invertible sheaf $L$ on $X$ is relatively semi-ample. 
More specifically, the following is our main result (note that it 
only holds  in positive characteristic, cf. \S \ref{s_c-zero}):

\begin{thm}\label{t_main}
Let $f\colon X \to S$ be a projective morphism 
of noetherian $\F_p$-schemes. 
Let $L$ be an invertible sheaf on $X$. 
Assume that  $L|_{X_s}$ is semi-ample for all the  
points  $s \in S$, where $X_s$ denotes the fibre of $f$ over $s$. 

Then $L$ is $f$-semi-ample. 
\end{thm}

In general, even if the schemes $X$ and $S$ appearing in Theorem \ref{t_main}, are of finite type over a field of positive characteristic, we need to consider not only closed points of $S$ 
but all the points of $S$ (cf. Example~\ref{e_scheme}). 
On the other hand, we may ignore non-closed points of $S$ 
if the base field is uncountable:

\begin{thm}\label{t_main2}
Let $k$ be an uncountable field of positive characteristic and 
let $f\colon X \to S$ be a projective $k$-morphism 
of schemes of finite type over $k$. 
Let $L$ be an invertible sheaf on $X$. 
Assume that $L|_{X_s}$ is semi-ample for all the closed points $s \in S$, where $X_s$ denotes the fibre of $f$ over $s$. 

Then $L$ is $f$-semi-ample. 
\end{thm}

\subsection{Description of the proof}

Although the schemes $X$ and $S$ appearing in Theorem \ref{t_main} could be of infinite dimension, it is easy to reduce the problem to the case where $X$ is of finite dimension (cf. Remark \ref{r-reduce-hensel}). 
Furthermore, replacing $S$ by $\Spec\,\widehat{\MO_{S, s}}$ for a point $s \in S$, 
we may assume that $X$ and $S$ are excellent. 
Then the proof of Theorem~\ref{t_main} proceeds by induction on the dimension of $X$. 
To clarify the structure of the proof, we introduce the following three statements:

\begin{theorema}\label{t-A}
Let $f\colon X \to S$ be a projective surjective morphism 
of excellent $\F_p$-schemes with connected fibres, where $X$ is normal
and of dimension $n \in \Z_{\geq 0}$. 
Let $L$ be an  invertible sheaf on $X$ such that 
$L|_{X_s}$ is semi-ample for all the  points $s \in S$. 

Then $L$ is $f$-semi-ample. 
\end{theorema}

\begin{theorema}\label{t-nume-triv4}
Let $f\colon X \to S$ be a projective surjective morphism of excellent reduced $\F_p$-schemes, where $X$ is of dimension $n \in \Z_{\geq 0}$. 
Let $L$ be an $f$-numerically trivial invertible sheaf on $X$ such that 
$L|_{X_s}$ is semi-ample for all the points $s \in S$. 

Then $L$ is $f$-semi-ample. 
\end{theorema}

\begin{theorema}\label{t-C}
Let $f\colon X \to S$ be a projective surjective morphism 
of excellent $\F_p$-schemes with connected fibres, where $X$ has dimension $n \in \Z_{\geq 0}$. 
Let $L$ be an  invertible sheaf on $X$ such that 
$L|_{X_s}$ is semi-ample for all the points $s \in S$. 

Then $L$ is $f$-semi-ample. 
\end{theorema}

\begin{rem}
After we submitted a preliminary version of this paper, B. Bhatt kindly informed us that he and P. Scholze have a  proof of Theorem \ref{t-nume-triv4} as a consequence of \cite[Theorem 1.3]{BS17}. Since their proof is very different from ours, we decided to keep it as it was (see Section \ref{s-nume-triv}). 
\end{rem}

For any $n \in \Z_{\geq 0}$, 
we denote by 
$({\rm Theorem~\ref{t-A}})_n$, 
$({\rm Theorem~\ref{t-nume-triv4}})_n$, or 
$({\rm Theorem~\ref{t-C}})_n$ 
the corresponding theorem in the case where $X$ has dimension $n$. 
For any $n, m \in \Z_{\geq 0}$, 
$({\rm Theorem~\ref{t-nume-triv4}})_{n, m}$
denotes the corresponding theorem in the case where $X$ has dimension $n$ 
and $S$ has dimension $m$. 
The proof of our main theorem is divided into three steps.

\medskip

\begin{enumerate}
\item[(I)]\label{sketch_I} 
(Theorem~\ref{t-C}$)_{n-1}$ implies (Theorem~\ref{t-A}$)_{n}$ (cf. Theorem~\ref{t_C-to-A}). 
\item[(II)]\label{sketch_II}   (Theorem~\ref{t-A}$)_{n}$ implies 
$({\rm Theorem~\ref{t-nume-triv4}})_n$ (cf. Theorem~\ref{t_A-to-B}). 
\item[(III)]\label{sketch_III} 
(Theorem~\ref{t-A}$)_n$ and $({\rm Theorem~\ref{t-nume-triv4}})_n$ imply (Theorem~\ref{t-C}$)_n$ (cf. Theorem~\ref{t_AB-to-C}). 
\end{enumerate}

\medskip
We now briefly describe these steps.
\medskip

(I)  
Let $f\colon X \to S$ be as in (Theorem~\ref{t-A})$_n$. 
As $X$ is normal, we may assume by standard arguments that 
both $X$ and $S$ are integral normal schemes. 
Using the Iitaka fibration induced by $L|_{X_{K(S)}}$ 
where $X_{K(S)}$ denotes the generic fibre of $f$, 
we are reduced to the case where $L|_{X_{K(S)}}$ is 
numerically trivial or ample (cf. Claim in the proof of Theorem \ref{t_C-to-A}). 
Note that, in this argument, we might replace $X$ by a birational model and this  
 requires the condition of  $X$ to be normal. 
If $L|_{X_{K(S)}}$ is numerically trivial, 
then we are done by taking a suitable alteration of the base scheme 
(cf. Proposition \ref{p-normal-nt1}). 
Thus, it suffices to treat the case where $L|_{X_{K(S)}}$ is ample. 
By a relative version of Keel's theorem (cf. Proposition~\ref{p-relative-keel}), it is enough to show that the restriction of $L$ to its $f$-exceptional locus $\mathbb E_f(L)$ is relatively semi-ample. 
This directly follows from (Theorem~\ref{t-C}$)_{n-1}$.

\medskip

(II) 
Let $f\colon X \to S$ be as in  (Theorem~\ref{t-nume-triv4})$_n$.
We may reduce the problem to the case where $S$ is 
an integral normal scheme (cf. Proposition \ref{p-reduction-normal}). 
Let $\nu\colon Y\to X$ be the normalisation of $X$, and let $C_X$ and $C_Y$ denote the conductors in $X$ and $Y$ respectively. 
Then we proceed by a quadruple induction on $(\dim X, \dim S, \delta(f), \eta(f)) 
\in \Z_{\geq 0}^4$, 
where we equip $\Z_{\geq 0}^4$ with the lexicographical order and, if $\overline{\xi}$ 
is the geometric generic point of $S$ and  $C_{X,\overline{\xi}}$ is the fibre of $C_X\to S$ over $\overline{\xi}$, 
 we denote by
$\delta(f)$ the dimension of $X_{\overline{\xi}}$   and by
$\eta(f)$   the number of the connected components of $C_{X, \overline{\xi}}$. 
As we are assuming  (Theorem~\ref{t-A}$)_{n}$, 
we have that $\nu^*L$ is relatively semi-ample and,
by the induction hypothesis, 
we may assume that $L|_{C_X}$ is relatively semi-ample.

By a result of Ferrand, we can normalise $X$ only
 along one horizontal component of $C_X$, 
which drops $\eta(f)$ exactly by one. 
For the sake of simplicity,  we briefly overview  two crucial cases: $\eta(f)=0$ and $\eta(f)=1$. 

Assume first that $\eta(f)=0$. 
After taking a suitable faithfully flat finite cover of $S$ (cf. Step \ref{step1-p-vertical} of Proposition \ref{p-vertical}), 
we may assume that there exists a closed subscheme $\Gamma$ of $X$ 
such that $\Gamma \to S$ is a generically universal homeomorphism. 
Applying Proposition \ref{p-thickening}, 
we may find a closed subscheme $X'$ on $X$ 
that is set-theoretically equal to $\Gamma$ over a generic locus over $S$ 
and which satisfies the following properties (cf. Step~\ref{step3-p-vertical} of Proposition~\ref{p-vertical}): 
\begin{enumerate}
\item[(i)] $L|_{X'}$ is relatively semi-ample by the induction hypothesis,  
and
\item[(ii)] the relative semi-ampleness of $L|_{X'}$ implies the one of $L$. 
\end{enumerate} 
Thus, we are done in the case $\eta(f)=0$. 

Assume now that  $\eta(f)=1$. We consider   the generic fibre $X_\eta$ of $f$ and, by assumption, the restriction of $L$ to $X_\eta$ is semi-ample. Using an argument similar to the previous case, we can show that $L$ is relatively semi-ample
(cf. Step~\ref{step5-p-nume-triv} of Theorem~\ref{t_A-to-B}). 
We refer to Section \ref{s-nume-triv}  for more details.

\medskip

(III) 
Let $f\colon X \to S$ be as in  (Theorem~\ref{t-C}$)_n$. 
We consider the normalisation $\nu\colon Y\to X$ of $X$. 
The most significant part of this case is to show that $L$ is EWM (cf. Subsection \ref{s_properties-invertible-sheaves}). 
To this end, inspired by \cite[Theorem 0.1]{keel03}, we prove  the following theorem (see Section \ref{s-AB-to-C} for its proof):

\begin{thm}\label{t-Keel-EWM}
Let $S$ be a noetherian $\F_p$-scheme. 
Let $f\colon Y \to X$ be a finite surjective 
$S$-morphism of reduced algebraic spaces proper over $S$.  
Let $L$ be an invertible sheaf on $X$ which is nef over $S$. 

Then $L$ is EWM over $S$ if and only if 
\begin{enumerate}
\item 
$L|_Y$ is EWM over $S$, and 
\item 
there exists a positive integer $m_0$ such that 
for all the geometric points $s \in S$, 
the $L|_{X_s}$-equivalence relation on $X_s$ is bounded by $m_0$ (cf. Definition~\ref{d_l-equivalence}).
\end{enumerate}
\end{thm}

By (Theorem~\ref{t-A})$_n$, 
we have that $f^*L=L|_Y$ is relatively semi-ample, 
hence (1) of Theorem \ref{t-Keel-EWM} holds. 
Moreover we have that (2) of Theorem \ref{t-Keel-EWM} also holds, 
by the assumption that $L|_{X_s}$ is semi-ample for all the  points $s \in S$.  
Therefore, we may apply Theorem \ref{t-Keel-EWM}, 
i.e. there exists an $S$-morphism $g\colon X\to Z$ to an algebraic space $Z$ proper over $S$ such that $g$ contracts all the $L$-trivial curves. 
By a variant of (Theorem \ref{t-nume-triv4})$_n$ (cf. Theorem \ref{t-nume-triv-as}), 
we conclude that $L^{\otimes m}=g^*L_Z$ 
for a positive integer $m$ and an invertible sheaf $L_Z$ on $Z$. 
Thus, the Nakai--Moishezon criterion implies that $Z$ is projective over $S$, as desired. 

\begin{rem}
It is worth explaining 
why the schemes which appear in Theorem~\ref{t-A}, \ref{t-nume-triv4}, and 
\ref{t-C}, 
are assumed to be not only noetherian but excellent. 
There are three advantages for this. 
First, it is necessary to impose the universally catenary condition 
to apply induction on the dimension of $X$ (cf. Section \ref{ss-catenary}). 
Second, we frequently take the normalisations of both the total and the base space, 
which compels us  to treat only universally Japanese schemes.  
Third, we  use Gabber's alteration theorem,  
which only holds for quasi-excellent schemes  
 (cf. Theorem \ref{t-alteration}). 
\end{rem}

\begin{rem}
Note that even if we are interested to prove Theorem \ref{t_main} only  for schemes of finite type over fields, 
our proof requires us to treat schemes that are not essentially of finite type over a field. 
This is because we repeatedly make use 
of henselian or complete local rings 
in the proof of ${\rm Theorem~\ref{t-nume-triv4}}$ (cf. Lemma \ref{l-finite/local}). 
\end{rem}

\section{Preliminary results}

\subsection{Notation and conventions}\label{s-nc} 

\begin{itemize}
\item 
A {\em variety} $X$ over a field $k$ is an integral scheme which is separated and of finite type over $k$. 
A {\em curve} is a variety of dimension one. Given a scheme $X$, we denote by $X_{\red}$ its {\em reduced structure}. We refer to \cite[I.\S 1]{hartshorne77} for the definition of dimension of a topological space. 
 
\item 
A morphism  $f\colon Y \to X$ of schemes is a {\em birational morphism} if there exists an open dense subset $X^0$ such that $f^{-1}(X^0)$ is dense in $Y$ and the induced morphism 
$f^{-1}(X^0) \to X^0$ is an isomorphism of schemes. 
\item 
Given a morphism $f\colon X \to Y$ of algebraic spaces, and given a point $y\in Y$, we denote by $X_y$ the fibre of $f$ over $y$. 
We say that $f$ has {\em connected fibres} or $f$ is a morphism 
with {\em connected fibres} if 
for any field $K$ and morphism $\Spec\,K \to Y$, 
the fibre product $X \times_Y \Spec\,K$ is a connected algebraic space. 
\item 
For definition of {\em catenary, universally catenary, quasi-excellent} and {\em excellent} schemes, 
we refer to \cite[Definition 8.2.1 and 8.2.35]{liu02}. 
Throughout this paper, 
excellent and quasi-excellent schemes are assumed to be quasi-compact i.e. noetherian, although  \cite[Definition 8.2.35]{liu02} does not impose such an assumption. 
\item 
An algebraic space $X$ is {\em noetherian} (resp. {\em excellent}) 
if $X$ is quasi-compact and for any \'etale morphism $U \to X$ 
from an affine scheme $U$, the ring $\Gamma(U, \MO_U)$ is 
a noetherian ring (resp. an excellent ring). 
Note that if $X \to Y$ is a morphism of finite type between algebraic spaces 
and  $Y$ is excellent, then so is $X$ (cf. \cite[\S32]{matsumura89} and \cite[Proposition 7.8.6]{egaIV.2}). 
\item 
Given an integral scheme $X$, we define $K(X):=\MO_{X, \xi}$ 
where $\xi$ is the generic point of $X$. 
For an integral domain $A$, we define $K(A):=K(\Spec\,A)$. 
\item 
Given an abelian group $H$, we define $H_{\Q}:=H \otimes_{\Z} \Q$ and given a homomorphism of abelian groups $\varphi\colon H\to K$, we denote by $\varphi_{\Q}\colon H_{\Q}\to K_{\Q}$ the induced homomorphism.
\item 
A morphism of noetherian schemes $f\colon X \to Y$ is
{\em generically finite} 
if there exists an open dense subset $Y'$ of $Y$ such that 
the induced morphism $f^{-1}(Y') \to Y'$ is a finite morphism 
(cf. \cite[Expos\'e II, Proposition 1.1.7 and the sentence after that]{ILO14}). 
\end{itemize}

\subsubsection{Properties of invertible sheaves}\label{s_properties-invertible-sheaves}

We refer to \cite{kollar13} for the classical definitions concerning a divisor on a proper normal varieties over a field $k$ (e.g. nef, semi-ample, big). 
Let $f\colon X\to S$ be a proper morphism of noetherian algebraic spaces 
and let $L$ be an invertible sheaf on $X$. 

\begin{itemize}
\item 
$L$ is $f$-{\em nef} if 
for any field $K$ and morphism $\Spec\,K \to Y$, 
the pullback of $L$ to the base change $X \times_Y \Spec\,K$ is nef 
(cf. Lemma~\ref{l-nef-special}). 
\item 
$L$ is $f$-{\em numerically trivial} if 
both $L$ and $L^{-1}$ are $f$-nef. 
\item 
$L$ is $f$-{\em free} if the natural homomorphism $f^*f_*L\to L$ is surjective.  
In particular, if $L$ is $f$-free then it induces a morphism $X\to \mathbb P(f_*L)$ over $S$. 
\item 
 $L$ is {\em $f$-very ample} if it is $f$-free and the induced morphism $X\to \mathbb P(f_*L)$ is a closed immersion. 
\item  
$L$ is $f$-{\em semi-ample} (resp. {\em  $f$-ample}) if $L^{\otimes m}$ is $f$-free (resp. $f$-very ample) for some positive integer $m$. 

\item $L$ is $f$-{\em weakly big} if there exist an $f$-ample invertible sheaf $A$ on $X$ 
and a positive integer $m$ such that if $g\colon X_{\red} \to S$ denotes the induced morphism, then  
\[
g_*((L^{\otimes m}\otimes_{\MO_X} A^{-1})|_{X_{\red}})\neq 0.
\]
Assuming that  $X$ is normal, $L$ is $f$-{\em big} if, for any connected component $Y$ of $X$, the restriction  $L|_{Y}$ is $h$-weakly big, where $h=f|_Y$ is the induced morphism. 
\item 
The $f$-{\em stable base locus} of $L$ is defined as  the following closed subset of $X$: 
\[
\mathbb B_f(L)=\bigcap_{m\ge 1} \Supp ~ {\rm Coker} (f^*f_*L^{\otimes m}\to L^{\otimes m}).
\]
\item 
Assume that $X$ is a scheme.
If $L$ is $f$-nef, the $f$-{\em exceptional locus} of $L$, denoted by $\mathbb E_f(L)$, is defined as the union of all the reduced closed subschemes  $V\subset X$ such that $L|_V$ is not $f|_V$-weakly big. 
Later, we shall  prove that $\mathbb E_f(L)$ is a closed subset of $X$ 
(cf. Lemma \ref{l_keel}). 
\item If $L$ is $f$-nef, then we say that $L$ is {\em endowed with a map} (EWM) over $S$ if there is a proper $S$-morphism $g\colon X\to Y$ to an algebraic space $Y$ proper 
over $S$ such that, for any point $s\in S$ and for any irreducible closed subspace $Z$ of $X_s$, we have that $\dim g(Z)<\dim Z$ if and only if $(L|_{X_s})^{\dim Z}\cdot Z=0$. 
\end{itemize}

When no confusion arises, if $L$ is $f$-nef (resp. $f$-big, \dots),  we will simply say that $L$ is relatively nef (resp. big, \dots) or $L$ is nef (resp. big, \dots) over $S$.

Note that if $X$ is a reduced scheme, then  \cite[Proposition 21.3.4]{egaiv.4} implies that any invertible sheaf on $X$  is of the form $\MO_X(D)$ where $D$ is a Cartier divisor on $X$.  

\subsubsection{Projective morphisms}

Let $f\colon X \to Y$ be a morphism of algebraic spaces. 
We refer to  \cite[Ch. II, Section 7]{knutson71}, for the definition of (quasi-)projective morphisms between algebraic spaces. 
If $X$ and $Y$ are schemes, these definitions coincide with the one in 
\cite[page 103]{hartshorne77}, but differ from the one given by Grothendieck 
\cite[D\'efinition 5.5.2]{egaii}. 
On the other hand, it is known that their definitions coincide in many cases 
(cf. \cite[Section 5.5.1]{fga2005}).

\subsection{Basic results}

In this subsection, we summarise some basic facts which will be used later. 
Although some of the material  here might be well-known, 
we provide their proofs for the sake of completeness. 

\begin{lemma}\label{l_connected-fibres}
Let $S$ be a noetherian $\F_p$-scheme and 
let $f\colon X \to Y$ be a surjective $S$-morphism of proper $S$-schemes with connected fibres. 

Then the induced map 
\[
H^0(Y,\MO_Y^{\times})_{\Q} \to H^0(X, \MO_{X}^{\times})_{\Q}
\]
is an isomorphism of groups.
\end{lemma}
\begin{proof}
Let 
\[
f\colon X\xrightarrow{f'} Y' \xrightarrow{\eta} Y
\]
be the Stein factorisation of $f$. Since the fibres of $f$ are connected,   
$\eta$ is a finite universal homeomorphism. 
By \cite[Proposition 6.6]{kollar97}, there exists a positive integer $e$ such that 
the $e$-th iterated Frobenius morphism $F^e\colon Y\to Y$ factors through $\eta$. 
Since $f'_*\MO_X=\MO_{Y'}$, it follows that  $H^0(Y',\MO_{Y'}^{\times})\to H^0(X, \MO_{X}^{\times})$ is bijective. Since $F^e$ factors through $\eta$, it follows that  $H^0(Y,\MO_Y^{\times})_{\Q} \to H^0(Y', \MO_{Y'}^{\times})_{\Q}$ is bijective. 
\end{proof}

\begin{lemma}\label{l_EWM}
Let $S$ be a noetherian $\mathbb F_p$-scheme and let $f\colon X\to Y$ be a finite universal homemorphism of algebraic spaces proper over $S$. Let $L$ be an invertible sheaf on $X$. 

Then $L$ is EWM over $S$ if and only if $f^*L$ is EWM over $S$. 
\end{lemma}
\begin{proof}
By \cite[Proposition 6.6]{kollar97}, there exists a positive integer $e$ such that 
the $e$-th iterated Frobenius morphism $F^e\colon X\to X$ factors through $f$. 
Thus, the claim follows.
\end{proof}

\begin{lemma}\label{l-basic-bc}
Let 
\[
\begin{CD}
X' @>\alpha>> X\\
@VVf'V @VVf V\\
S' @>\beta>> S
\end{CD}
\]
be a cartesian diagram of morphisms of schemes, where $\beta$ is an affine morphism. 

Then the induced homomorphism $$\theta\colon f^*\beta_*\MO_{S'} \to \alpha_*\MO_{X'}$$ is an isomorphism. 
\end{lemma}

\begin{proof}
We may assume that $S$ and $S'$ are affine: $S=\Spec\,R$, $S'=\Spec\,R'$. 
If $j\colon U\to X$ is an open immersion and $U':=U \times_X X'$, 
then we obtain 
\[
j^*\theta\colon j^*f^*\beta_*\MO_{S'} \to j^*\alpha_*\MO_{X'}=(\alpha|_{U'})_*\MO_{U'}.
\] 
Thus, we may assume that $X$ is affine: $X=\Spec\,A$. 
Then both sides of  
\[
\theta(X)\colon \Gamma(X, f^*\beta_*\MO_{S'}) \to 
\Gamma(X, \alpha_*\MO_{X'})
\]
are naturally isomorphic to $A \otimes_R R'$. 
Therefore $\theta$ is an isomorphism. 
\end{proof}

\begin{lemma}\label{l-local-flattening}
Let $A \subset B$ be an integral extension of integral domains 
such that the induced field extension $K(A) \subset K(B)$ is a finite extension. 

Then there exists a subring $B'$ of $K(B)$ which satisfies the following properties:  
\begin{enumerate}
\item $K(B')=K(B)$, and  
\item $B'$ contains $A$ and $B'$ is a free $A$-module whose rank is equal to $[K(B):K(A)]$. 
\end{enumerate}
\end{lemma}

\begin{proof}
We may assume that $K(A) \subset K(B)$ is a simple extension. 
Since $K(A) \subset K(B)$ is simple, 
there exists an element $\beta \in K(B)$ such that $K(B)=K(A)[\beta]$ and  
\[
\beta^n+\alpha_1\beta^{n-1}+\cdots+\alpha_n=0
\]
where $n:=[K(B):K(A)]$ and  $\alpha_1,\dots,\alpha_n \in K(A)$. 
For each $i$, we may write $\alpha_i=a_i/a'_i$ for some $a_i, a'_i \in A$ with $a'_i \neq 0$. 
Killing the denominators and after possibly replacing $\beta$ by $a\beta$ for some $a\in A \setminus \{0\}$, we may assume that $\alpha_i \in A$ for all $i$. 
In particular, $\beta$ is an element of $K(B)$ which is integral over $A$. 
Let
\[
B':=A[\beta].
\]
Consider the surjective $A$-algebra homomorphism
\[
\varphi\colon A[t] \to A[\beta]=B'\quad \text{such that }\varphi(t)=\beta.
\]
It is enough to show that $\Ker(\varphi)=f(t)A[t]$, where  
\[
f(t):=t^n+\alpha_1t^{n-1}+\cdots+\alpha_n \in A[t].
\]
Since the inclusion $\Ker(\varphi)\supset f(t)A[t]$ is obvious, 
it is enough to prove that  $\Ker(\varphi)\subset f(t)A[t]$. 
Pick $g(t) \in \Ker(\varphi)$. 
Since $f(t)$ is monic, we have 
\[
g(t)=f(t)h(t)+\sum_{i=0}^{n-1}c_i t^i
\]
for some $h(t) \in A[t]$ and $c_0,\dots,c_{n-1} \in A$. 
It follows that 
\[
0=g(\beta)=\sum_{i=0}^{n-1}c_i \beta^i.
\]
Since $1, \beta, \cdots, \beta^{n-1}$ is a $K(A)$-linear basis of $K(B)$, 
we obtain $c_0=c_1=\cdots=c_{n-1}=0$ and $g(t) \in f(t)A[t]$, as desired. 
\end{proof}

\begin{lemma}
\label{l-exc-fin-sub}
Let $R$ be a noetherian ring and let 
$A \subset B$ be a ring extension of $R$-algebras, 
where $B$ is a finitely generated $A$-module and a finitely generated $R$-algebra.  

Then  $A$ is a finitely generated $R$-algebra. 
\end{lemma}

\begin{proof}
Let $b_1, \cdots, b_m$ be generators of $B$ as an $R$-algebra. 
Since $A \subset B$ is an integral extension, for any $i \in \{1,\dots,m\}$, 
there exist $a_{i, 1},\dots,a_{i,n_i} \in A$ such that 
\[
b_i^{n_i}+a_{i, 1} b_i^{n_i-1}+...+a_{i, n_i}=0.
\]
Let  $A'$  be the $R$-subalgebra of $A$ generated by all the 
$a_{i, j}$. 
In particular, $A'$ is a finitely generated $R$-algebra.
 We have the  inclusions: 
\[
A' \subset A \subset B.
\]
Since $A'$ is a noetherian ring and 
$B$ is a finitely generated $A'$-module, 
also $A$ is a finitely generated $A'$-module.
Thus, $A$ is a finitely generated $R$-algebra, as desired. 
\end{proof}

\begin{lemma}\label{l-nef-special}
Let $f\colon X \to S$ be a proper morphism of noetherian schemes 
and let $L$ be an invertible sheaf on $X$. 

Then the following are equivalent:
\begin{enumerate}
\item $L$ is $f$-nef. 
\item $L|_{X_s}$ is nef for all the  points $s\in S$. 
\item $L|_{X_s}$ is nef for all the closed points $s\in S$. 
\end{enumerate}
\end{lemma}

\begin{proof}
It is enough to show that (3) implies (2). 
To this end, we may assume that $S=\Spec\,R$ 
where $R$ is a discrete valuation ring. 
Moreover, by Chow's lemma, 
we may assume that $f$ is projective. 

Let $\xi \in S$ (resp. $0 \in S$) be the non-closed (resp. closed) point. 
Given a curve $C_{\xi}$ on $X_{\xi}$ which is projective over $k(\xi)$, 
it is enough to show that $(L|_{X_{\xi}}) \cdot C_{\xi} \geq 0$. 
Since $f$ is projective, there exists a closed immersion $C \to X$ 
such that the composite morphism $C \to X \to S$ is flat and 
$C \times_S \Spec\,k(\xi)=C_{\xi}$. 
Since the intersection number is invariant under flat family, 
we get 
$$(L|_{X_{\xi}}) \cdot C_{\xi}=(L|_{X_0}) \cdot (C|_{X_0}) \geq 0,$$
as desired. 
\end{proof}

\subsection{Dimension formulas for universally catenary schemes }\label{ss-catenary}

The goal of this subsection is to show that some of the standard dimension formulas for a proper morphism between  varieties extend to the category of universally catenary schemes. 

We believe that the results in this subsection are well known, but we include proofs for completeness. 

\begin{lemma}\label{l-cate-dim}
Let $f\colon X \to Y$ be a proper surjective morphism 
of universally catenary noetherian integral schemes.

Then 
\[
\dim X=\dim Y+{\rm tr.deg}_{K(Y)}\,K(X).
\]
\end{lemma}

\begin{proof}
See \cite[Corollaire 5.6.5]{egaIV.2}. 
\end{proof}

\begin{prop}\label{p-dim-codim}
Let  $f\colon X \to Y$ be a proper surjective morphism 
of universally catenary noetherian integral schemes, where 
 $A$ is a local ring  and $Y=\Spec\,A$.  Let $X'$ be an irreducible closed subset of $X$. 

Then  there exists a sequence of irreducible closed subsets 
of $X$
\[
X=:X_{\dim X} \supsetneq X_{\dim X-1} \supsetneq \cdots \supsetneq X_0  \neq \emptyset
\]
such that $X'=X_i$ for some $i \in \{0, \cdots, \dim X\}$. 

In particular, 
\[
\dim X'+\codim_X X'=\dim X.
\]
\end{prop}

\begin{proof}
We first treat the case where  $X'=\{x\}$ 
 for some closed point $x$ of $X$. 
Since $f$ is proper, 
the image $y:=f(x)$ is a closed point of $Y$. 
Then we have that 
\[
\dim \MO_{X, x}-\dim \MO_{Y, y}={\rm tr.deg}_{K(Y)} K(X)
=\dim X-\dim Y
\]
where the first (resp. the second) equality  
holds by \cite[Theorem 8.2.5]{liu02} 
(resp. Lemma~\ref{l-cate-dim}). 
As $\codim_X \{x\}=\dim \MO_{X, x}$ and 
$\dim \MO_{Y, y}=\codim_Y \{y\}=\dim Y$, the claim follows. 

We now prove the general case. 
We fix a closed point $x$ of $X$ which is contained in $X'$. 
Then $X'$ corresponds to a prime ideal $\p$ of 
the local ring $\MO_{X, x}$ at $x$. 
Since the claim holds in the case $X'=\{x\}$,  we have that 
$\dim \MO_{X, x}=\dim X$. 
Thus, 
\[
\dim (\MO_{X, x}/\p)+\dim (\MO_{X, x})_{\p}=\dim \MO_{X, x}=\dim X,
\]
where the first equality follows from 
the fact that $\MO_{X, x}$ is catenary. Thus, the claim follows. 
\end{proof}

Below, given a  morphism $f\colon X\to Y$ between schemes and given a subset $W$ of $X$ (resp. $W'$ of $Y$) we denote by $f(W)$ (resp. $f^{-1}(W')$) the set-theoretic image (resp. inverse image) of $W$ (resp. $W'$). 

\begin{lemma}\label{l-prime-push}
Let $f\colon X \to Y$ be a proper surjective morphism 
of noetherian universally catenary schemes. 
Let $r:=\dim X-\dim Y$. 
Assume that $f^{-1}(y)$ is pure $r$-dimensional for any closed point $y \in Y$. 

Then the following hold:
\begin{enumerate}
\item 
For any irreducible closed subset $Y_1$ of $Y$ and 
any irreducible component $X_1$ of $f^{-1}(Y_1)$ satisfying $f(X_1)=Y_1$, 
we have that $\dim X_1-\dim Y_1=r.$ 
\item 
Assume that $X$ and $Y$ are integral schemes. 
If $D$ is an irreducible closed subset such that $\codim_X D=1$, 
then 
\[
\codim_Y f(D) \leq 1.
\]
\end{enumerate}
\end{lemma}

\begin{proof}
We first show (1). Let $Y_1$ and $X_1$ be as in the statement. 
We may assume that $\dim Y_1<\infty$ and we prove the claim by induction on $\dim Y_1$. 
If $\dim Y_1=0$, then there is nothing to show. 
Thus, we may assume that $\dim Y_1>0$.
By  generic flatness, there exists a point $z \in Y_1$ such that  
\[
\dim X_1-\dim Y_1=\dim (f^{-1}(z) \cap X_1) \leq \dim f^{-1}(z)=r.
\]
Thus, it is enough to show that $\dim X_1-\dim Y_1 \geq r$. 
As $\dim Y_1>0$, 
we can find an irreducible closed subset $Y_2$ of $Y_1$ satisfying 
$\dim Y_2=\dim Y_1-1$. 
Since
\[
f(X_1 \cap f^{-1}(Y_2))=f(X_1) \cap Y_2=Y_2,
\]
there is an irreducible component $X_2$ of $X_1 \cap f^{-1}(Y_2)$ 
such that $f(X_2)=Y_2$. 
By  induction, it follows that $\dim X_2-\dim Y_2=r$. 
Since 
\[
X_2 \subset X_1 \cap f^{-1}(Y_2) \subsetneq X_1,
\]
we have that $\dim X_2 <\dim X_1$. Thus, 
\[
\dim X_1-\dim Y_1 \geq (\dim X_2+1)-(\dim Y_2+1)=r
\]
and (1) holds. 

We now show (2). Let $y$ be the generic point of $f(D)$. After 
replacing $f\colon X \to Y$ by the base change 
$X \times_Y \Spec\,\MO_{Y, y} \to \Spec\, \MO_{Y, y}$ 
we may assume that $Y=\Spec\,A$ for some local ring $A$. 
If $f(D)=Y$, then there is nothing to show. 
Thus, we may assume that $f(D) \subsetneq Y$. 
By (1), we have that  $\dim f^{-1}(f(D))<\dim X$.  Since $\codim_X D=1$, it follows that  $D$ is an irreducible component of $f^{-1}(f(D))$. Proposition~\ref{p-dim-codim} implies 
\[
\codim_Y f(D)=\dim Y-\dim f(D),
\]
and
\[ 1=\codim_X D=\dim X-\dim D.
\]
Since $D$ is an irreducible component of $f^{-1}(f(D))$, we have 
\[
\codim_Y f(D)=\dim Y-\dim f(D)=(\dim X-r)-(\dim D-r)=1,
\]
where the second equality follows from (1). 
Thus, (2) holds. 
\end{proof}

\subsection{Relative semi-ampleness}
The purpose of this subsection is to recall some basic results on the relative semi-ampleness of an invertible sheaf. Many of these results are well-known however we provide proofs for the sake of completeness. 

\begin{lemma}\label{l-stein}
Let 
\[
f\colon X \overset{f'}\to S' \overset{\alpha}\to S
\] 
be proper morphisms of noetherian schemes and 
let $L$ be an invertible sheaf on $X$.

Then the following hold:
\begin{enumerate}
\item If $L$ is $f$-semi-ample, then  $L$ is $f'$-semi-ample. 
\item If $L$ is $f'$-semi-ample and $\alpha$ is finite, then $L$ is $f$-semi-ample. 
\end{enumerate}
\end{lemma}

\begin{proof}
For any positive integer $m$, we have 
\[
f^*f_*L^{\otimes m}=f'^*\alpha^*\alpha_*f'_*L^{\otimes m}  \to f'^*f'_*L^{\otimes m}\to L^{\otimes m}.
\]
Thus, (1) holds. 
Since $\alpha$ is finite, we have that 
\[
f^*f_*L^{\otimes m} =f'^*\alpha^*\alpha_*f'_*L^{\otimes m}  \to f'^*f'_*L^{\otimes m} 
\]
is surjective.  Thus, (2) holds. 
\end{proof}

\begin{lemma}\label{l_pullback}
Let 
\[
f'\colon X' \overset{\beta }\to X \overset{f}\to S
\] 
be proper morphisms of noetherian schemes and let $L$ be an invertible sheaf on $X$. 

Then the following hold:
\begin{enumerate}
\item 
If $L$ is $f$-semi-ample, then $\beta^*L$ is $f'$-semi-ample. 
\item 
If $\beta_*\MO_{X'}=\MO_X$ and $\beta^*L$ is $f'$-semi-ample, then $L$ is $f$-semi-ample. 
\item 
If $X$ is an $\F_p$-scheme, $\beta$ has connected fibres and $\beta^*L$ is $f'$-semi-ample, then $L$ is $f$-semi-ample. 
\item If $S$ is excellent,
$X$ is normal, $\beta$ is surjective and $\beta^*L$ is $f'$-semi-ample, then $L$ is $f$-semi-ample. 
\end{enumerate}
\end{lemma}
\begin{proof}
If $L$ is $f$-semi-ample, there is a positive integer $m$ such that 
\[
f^*f_*L^{\otimes m}\to L^{\otimes m}
\]
is surjective. 
Thus, the composite morphism 
\[
\beta^*f^*f_*L^{\otimes m}=f'^*f_*L^{\otimes m}\to f'^*f_*\beta_*\beta^*L^{\otimes m}= f'^*f'_*\beta^*L^{\otimes m} \to \beta^* L^{\otimes m}
\]
is surjective. In particular, $f'^*f'_*\beta^*L^{\otimes m} \to \beta^* L^{\otimes m}$ is surjective.
Thus, (1) holds. 

We now show (2). 
To this end, we may assume that $S$ is affine. 
Pick a closed point  $x\in X$. Then, since $\beta$ is proper and surjective,  there exist a closed point $x'\in X'$, a positive integer $m$  and $t\in H^0(X',\beta^*L^{\otimes m})$ such that $\beta|_{x'}=x$ and $t|_{x'}\neq 0$. Since $\beta_*\MO_{X'}=\MO_X$, there exists $s\in H^0(X , L^{\otimes m})$ such that $s|_x\neq 0$. It follows that $L$ is semi-ample over $S$. 
Thus, (2) holds.

We now show (3). 
Let $X' \to X'' \to X$ be the Stein factorisation of $\beta$. Since the fibres of $\beta$ are connected, 
we have that $X'' \to X$ is a universal homeomorphism. 
Thus, by (2), we may assume that $\beta$ is a universal homeomorphism. 
By \cite[Proposition 6.6]{kollar97}, there exists a positive integer $e$ such that 
the $e$-th iterated Frobenius morphism $F^e\colon X\to X$ factors through $\beta$. 
Hence, replacing $\beta$ by $F^e$, we may assume that $\beta=F^e$. 
In this case, the assertion (3) is clear.

We now show (4). 
Taking the Stein factorisation of $\beta$, 
(2) implies that we may assume that $\beta$ is a finite morphism.  
Moreover, replacing $X'$ by its normalisation, 
the problem is reduced to the case where $X'$ is normal. 
If the field extension $K(X) \subset K(X')$ is purely inseparable, 
then the assertion follows from (3). 
Therefore, taking the separable closure of $K(X) \subset K(X')$, 
we see that the problem is reduced to the case where 
the field extension of $K(X) \subset K(X')$ is separable. 
Furthermore, taking its Galois closure, we may assume that 
$K(X) \subset K(X')$ is a Galois extension with Galois group $G$. 
Pick a closed point  $x\in X$ and let $\{x'_1,\dots,x'_k\}$ be the inverse image of $x$ by $\beta$. There exist a positive integer $m$ and $t\in H^0(X',\beta^*L^{\otimes m})$ such that $t(x'_i)\neq 0$ for any $i\in \{1,\dots,k\}$. 
Then 
\[
t':=\prod_{\sigma \in G} \sigma^*t_i\in H^0(X',\beta^*L^{\otimes m|G|})
\]
descends to $X$, i.e. there exists $s\in H^0(X,L^{\otimes m|G|})$ such that $\beta^*s=t$. 
In particular, $s|_x \neq 0$. Thus (4) holds. 
\end{proof}

\begin{lemma}\label{l-ff-descent}
Let 
$$\begin{CD}
X' @>\beta >> X\\
@VVf'V @VVfV\\
S' @>\alpha >> S
\end{CD}$$
be a cartesian diagram of morphisms of noetherian schemes, 
where $f$ is proper. 
Let $L$ be an invertible sheaf on $X$ and let $L':=\beta^*L$. 

Then the following hold:
\begin{enumerate}
\item 
If $L$ is $f$-semi-ample, then $L'$ is $f'$-semi-ample. 
\item 
If $L'$ is $f'$-semi-ample and $\alpha$ is faithfully flat, 
then $L$ is $f$-semi-ample. 
\end{enumerate}
\end{lemma}

\begin{proof} 
By (1) of Lemma~\ref{l_pullback}, if $L$ if $f$-semi-ample, then $L'$ is $(f\circ \beta)$-semi-ample. By (1) of Lemma~\ref{l-stein}, we have  that  $L'$ is $f'$-semi-ample. 
Thus, (1) holds. 

We now show (2). 
Since $L'$ is $f'$-semi-ample, 
there exists a positive integer $m$ such that 
$f'^*f'_*L'^{\otimes m} \to L'^{\otimes m}$ is surjective. 
Since $\beta$ is faithfully flat, 
it suffices to show that $\beta^*(f^*f_*L^{\otimes m}) \simeq f'^*f'_*L'^{\otimes m}$, 
which follows from \cite[Proposition III.9.3]{hartshorne77}. 
Thus, (2) holds. 
\end{proof}

\begin{lemma}\label{l_reduced}
Let $f\colon X\to S$ be a proper morphism of noetherian $\mathbb F_p$-schemes and let $L$ be an invertible sheaf on $X$. 
Let $f':X_{\red} \xrightarrow{j} X \xrightarrow{f} S$, 
where $j$ is the induced closed immersion. 

Then 
\[
\mathbb B_f(L)=\mathbb B_{f'}(L|_{X_{\red}}).
\]
In particular, $L$ is $f$-semi-ample if and only if $L|_{X_{\red}}$ is $f'$-semi-ample. 
\end{lemma}

\begin{proof}
We may assume that $S$ is affine. 
Clearly, $\mathbb B_{f'}(L|_{X_{\red}})\subset \mathbb B_f(L)$. We now show the opposite inclusion. 
 Let $x\in X$ be a closed point such that $x\notin \mathbb B_{f'}(L|_{X_{\red}})$. Then there exist a positive integer $m$ and $s\in H^0(X_{\red},L^{\otimes m}|_{X_{\red}})$ such that $s|_x \neq 0$. Let  $F\colon X\to X$ be the absolute Frobenius morphism.  There exists a positive integer $e$ such that 
if $t=(F^e)^*(s)$, then $t\in H^0(X,L^{\otimes mp^e})$ and $t|_x \neq 0$. Thus, the claim follows. 
\end{proof}

\begin{rem}\label{r-reduce-hensel}
Let $f\colon X \to S$ be a proper morphism of noetherian schemes. 
Let $L$ be an invertible sheaf on $X$. 
Then the following are equivalent:
\begin{enumerate}
\item 
$L$ is $f$-semi-ample. 
\item 
For any point $s \in S$, 
if $S':=\Spec\,\MO_{S, s} \to S$ is the induced morphism and 
$\alpha\colon X':=X \times_S S' \to X$ is the projection, 
then 
$\alpha^*L$ is semi-ample over $S'$. 
\item 
For any point $s \in S$, 
if $S'':=\Spec\,\MO^h_{S, s} \to S$ is the induced morphism 
for the henselisation $\MO_{S, s}^h$ 
and 
$\beta\colon X'':=X \times_S S'' \to X$ is the projection, 
then 
$\beta^*L$ is semi-ample over $S''$. 
\item 
For any point $s \in S$, 
if $S''':=\Spec\,\widehat{\MO_{S, s}} \to S$ is the induced morphism 
for the completion $\widehat{\MO_{S, s}}$ 
and 
$\gamma:X''':=X \times_S S''' \to X$ is the projection, 
then 
$\gamma^*L$ is semi-ample over $S'''$. 
\end{enumerate}
Indeed, it is clear that (1) and (2) are equivalent. 
It follows from Lemma~\ref{l-ff-descent} that (2), (3) and (4) are equivalent. 
\end{rem}

\begin{lemma}\label{l-descend} 
Let $f\colon X \to S$ be a proper surjective morphism of 
noetherian $\F_p$-schemes with connected fibres. 
Let $L$ be an invertible sheaf which is 
$f$-numerically trivial and $f$-semi-ample. 

Then there exists a positive integer $m$ and 
an invertible sheaf $M$ on $S$ such that 
$L^{\otimes m} \simeq f^*M$. 
\end{lemma}

\begin{proof}
We can apply the same proof as in \cite[Lemma 1.1]{keel99}. 
\end{proof}

\begin{lemma}\label{l-finite/local}
Let $f\colon X \to S=\Spec\,R$ be a proper morphism of 
noetherian schemes 
and assume that there is a finite ring homomorphism $R_0 \to R$ 
such that $R_0$ is a henselian local ring. 
Let $L$ be  an invertible sheaf  on $X$. 

Then  the following are equivalent:
\begin{enumerate}
\item 
There exists a positive integer $m$ such that $L^{\otimes m} \simeq \MO_X$. 
\item 
$L$ is $f$-semi-ample and $f$-numerically trivial. 
\end{enumerate}
\end{lemma}

\begin{proof}
It suffices to show that (2) implies (1). 
Let 
$f\colon X \overset{g}\to T \to S$
 be the Stein factorisation of $f$. 
By Lemma \ref{l-descend}, there exists  a positive integer $m$ such that $L^{\otimes m} \simeq g^*M$ 
for some invertible sheaf $M$ on $T$. 
We can write $T=\Spec\,A$ for some ring $A$ finite over $R$, hence also over $R_0$. 
By \cite[Proposition 2.8.3]{Fu15}, $A$ is the direct product of finitely many local rings. 
Thus,  $M$  is trivial, and in particular 
also $L^{\otimes m}$ is trivial. 
\end{proof}

For notational convenience, we state the lemma below using Cartier divisors instead of invertible sheaves.

\begin{lemma}\label{l-flat-descent}
Let $f\colon X \to S$ be a proper morphism 
of integral normal excellent schemes satisfying $f_*\MO_X=\MO_S$.
Let $L$ be a $\Q$-Cartier $\Q$-divisor on $X$. 
Assume that 
\begin{enumerate}
\item{$S$ is $\Q$-factorial.}
\item{$L$ is $f$-nef.}
\item{$L|_{X_{K(S)}}\sim_{\Q}0$.}
\item For any prime divisor $D$ on $X$, 
its image $f(D)$ is either equal to $S$ or a prime divisor on $S$. 
\end{enumerate}

Then there exists a $\Q$-Cartier $\Q$-divisor $M$ on $Y$ such that $L\sim_{\Q} f^*M$. 
\end{lemma}

\begin{proof}
After possibly replacing $L$ by $rL$ for some positive integer $r$, 
we may assume that $L$ is a Cartier divisor. 
By (3),  we may find  a positive integer $m$ and 
$\varphi \in K(X)$ such that 
\[
mL+{\rm div}(\varphi)=L',
\]
where $L'$ is a Cartier divisor on $X$ such that 
$\Supp\,L' \subset f^{-1}(S^0)$ for some proper closed subset $S^0$ of $S$.

We show the claim by induction on the number of  irreducible components of $f(L')$. 
If this number is zero i.e. if  $L'=0$, then there is nothing to show. 
Thus, we may assume that $L' \neq 0$. Let $D$ be a prime divisor which is contained in the support of $L'$. Let $E:=f(D)$. Then (4) implies that  $E$ is a prime divisor and  (1) implies that $E$ is $\Q$-Cartier. 
We may write 
\[
f^*E=\sum_{i \in I} e_iD_i,
\]
where, for each $i\in I$, $D_i$ is a prime divisor and $e_i$ is a positive rational number. 
There exists a unique rational number $\alpha \in \Q$ such that if 
\[
L'':=L'-\alpha f^*E,
\]
then the coefficient of $L''$ along  $D_i$ is non-positive for any $i\in I$ and the coefficient of $L''$ along $D_{i_1}$ is equal to zero for some $i_1 \in I$. 
We define
\[
I':= \{i \in I\,|\, \text{the coefficient of } L''  \text{ along } D_i \text{ is negative}\}.
\]

We distinguish two cases. We first assume that
$I'=\emptyset$. Then
the number of  irreducible components of $f(L'')$ 
is less than the one of $f(L')$. 
By induction, it follows that  $L \sim_{\Q} f^*M$ for some $M$. 
Thus, we are done.

We now assume that $I' \neq \emptyset$. 
We want to derive a contradiction. 
By (4), for each $i\in I$, we have that $D_i$ dominates $E$. 
Let $K:=K(E)$. 
By abuse of notation, we denote  by $K$ also the generic point of $E$. 
The fibre $X_K$ of $X \to S$ over $K$ may be written as 
\[
X_K=\bigcup_{i \in I}\Supp\,(D_i)_K
=\left(\bigcup_{i \in I\setminus I'}\Supp\,(D_i)_K\right) \cup 
\left(\bigcup_{i \in I'}\Supp\,(D_i)_K\right).
\]
Since $X_K$ is connected, 
we can find $j_1 \in I\setminus I'$ and $j_2 \in I'$ such that 
$(D_{j_1})_K \cap (D_{j_2})_K \neq\emptyset$.

Since the coefficient of $-L''$ along any prime divisor intersecting $X_K$ 
is non-negative, 
there exists an open neighbourhood $\widetilde S$ of $K \in S$ 
such that $-L''|_{\widetilde X}$ is effective, where $\widetilde X:=f^{-1}(\widetilde S)$. 
Fix a positive integer $\ell$ such that $\ell L''$ is a Cartier divisor. 
Let 
\[
s \in H^0(\widetilde{X}, \MO_{\widetilde X}(-\ell L''))
\]
be the section corresponding to the effective Cartier divisor $-\ell L''|_{\tilde X}$. 
In particular,  $s|_{D_{j_1} \cap \widetilde{X}} \neq 0$ and $s|_{D_{j_2} \cap \widetilde{X}}= 0$. 
Thus,  $s|_{(D_{j_1})_K} \neq 0$ and $s_{(D_{j_2})_K} = 0$. 
Since $(D_{j_1})_K \cap (D_{j_2})_K \neq\emptyset$, 
we can find a $K$-curve $C$
 such that  
$s|_C \neq 0$ and $C \cap D_{j_2} \neq \emptyset$. 
In particular,
\[
s|_C \in H^0(C, \MO_C(-\ell L''))
\]
is such that $s|_z =0$ for any point $z \in C \cap D_{j_2}$. 
Thus, $\deg_C(-L''|_C)>0$, and in particular 
$L\cdot C= L''|_{X_K} \cdot C<0$, 
which contradicts the assumption that $L$ is $f$-nef. 
\end{proof}

\subsection{Relative Keel's theorem}

The goal of this subsection is to prove a relative version of Keel's theorem \cite[Theorem 0.2]{keel99}. 
To this end, we follow similar methods as in \cite[Lemma~3.3]{cmm11}. 

We begin with the following: 

\begin{lemma}\label{l_keel}
Let $f\colon X\to S$ be a projective surjective morphism of noetherian $\F_p$-schemes. Let $L$ be a $f$-nef invertible sheaf on $X$. 

Then the following hold: 
\begin{enumerate}
\item 
Given an $f$-ample invertible sheaf $A$, 
a positive integer $m$ and an element 
$s \in  H^0(X_{\red}, (L^{\otimes m}\otimes_{\MO_X} A^{-1})|_{X_{\red}})$,
if $Z$ is the reduced closed subscheme of $X$ 
whose support is equal to the zero set of $s$, and  
$g\colon Z \hookrightarrow X \xrightarrow{f} S$ is the induced moprhism, then   $\mathbb E_f(L)=\mathbb E_g(L|_Z)$.
\item $\mathbb E_f(L)=X$ if and only if $L$ is not $f$-weakly big.
\item  $\mathbb E_f(L)$ is a closed subset of $X$. 
\end{enumerate}
\end{lemma}

\begin{proof}
We first show (1) and (2). 
Clearly, the inclusion $\mathbb E_f(L) \supset \mathbb E_g(L|_Z)$ holds. 
Thus, it is enough to show the opposite inclusion. 
Pick a reduced closed subscheme $V$ of $X$ such that 
$L|_V$ is not $f|_V$-weakly big. 
Then $s|_V\in H^0(V, L|_V^{\otimes m}\otimes A^{-1}|_V)$ is equal to zero. 
It follows that $\Supp\, V\subset \Supp\, Z$, which implies that 
$\mathbb E_f(L) \subset \mathbb E_g(L|_Z)$. Thus, (1) holds.

Note that if $U\subset S$ is an open subset and if $f'\colon X':=X\times_S U\to U$ is the projection, then $\mathbb E_{f'}(L|_{X'})= \mathbb E_f(L)\cap X'$. 
Thus, in order to prove (2) and (3), we may assume that $S$ is affine. In this case,  (2)  follows immediately from (1). 

We now show (3). 
By (2), we may assume that $L$ is $f$-weakly big. 
Thus, there exist an $f$-ample invertible sheaf $A$, 
a positive integer $m$ and a nonzero element 
$s \in  H^0(X_{\red}, (L^{\otimes m}\otimes_{\MO_X} A^{-1})|_{X_{\red}}).$ 
Let $Z$ and $g$ be as in (1). 
Then $Z$ is a closed subscheme of $X$ such that  
$\Supp ~Z\subsetneq \Supp ~X$ and (1) implies that $\mathbb E_f(L)=\mathbb E_g(L|_Z)$. 
By noetherian induction, 
we may assume that $\mathbb E_g(L|_Z)$ is a closed subset of $Z$. 
Hence it is also a closed subset of $X$. 
Thus, (3) holds. 
\end{proof}

\begin{lemma}\label{l-relative-keel}
Let $f\colon X \to S$ be a projective morphism of noetherian 
$\F_p$-schemes, where $S$ is affine. 
Let $L$ be an $f$-nef invertible sheaf  on $X$ and 
let $D$ be an effective Cartier divisor on $X$ such that 
$A:=L \otimes_{\MO_X} \MO_X(-D)$ is $f$-ample. 
Let $r$ be a positive integer and let $t \in H^0(D, L^{\otimes r}|_D)$. 

Then there exists a positive integer $e_0$ and $t'\in H^0(X, L^{\otimes rp^{e_0}})$ such that 
$t'|_{p^{e_0}D}=(F^{e_0})^*t$, where $F^{e_0}\colon p^{e_0}D\to D$ 
is the morphism induced by the $e_0$-th iterated 
absolute Frobenius morphism $F^{e_0}\colon X \to X$. In particular, $t'|_D=t^{\otimes p^{e_0}}$.
\end{lemma}

\begin{proof}
Consider the exact sequence
\[
0\to L^{\otimes r}\otimes_{\MO_X}\MO_X(-D) \to L^{\otimes r}\to L^{\otimes r}|_D\to 0.
\]
For any positive integer $e$, we obtain the exact sequence
\[
0\to L^{\otimes rp^e}\otimes_{\MO_X}\MO_X(-p^eD) \to L^{\otimes rp^e}\to L^{\otimes rp^e}|_{p^eD}\to 0
\]
induced by taking the pull-back by the Frobenius morphism $F^e\colon X\to X$. 
Since $L$ is $f$-nef and $A$ is $f$-ample, it follows that the invertible sheaf 
\[
L^{\otimes r}\otimes_{\MO_X}\MO_X(-D)\simeq L^{\otimes (r-1)} \otimes_{\MO_X} A
\]
is $f$-ample. 
In particular, we can find a positive integer $e_0$ such that  
\[
H^1(X, L^{\otimes rp^{e_0}}\otimes_{\MO_X}\MO_X(-p^{e_0}D))\simeq 
H^1(X, (L^{\otimes r}\otimes_{\MO_X}\MO_X(-D))^{\otimes p^{e_0}})=0.
\]
Thus,
\[
H^0(X,L^{\otimes rp^{e_0}})\to H^0(X,L^{\otimes rp^{e_0}}|_{p^{e_0}D})
\]
is surjective. 
Therefore, 
there exists $t' \in H^0(X, L^{\otimes rp^{e_0}})$ such that 
$t'|_{p^{e_0}D}=(F^{e_0})^*t$, as claimed. 
\end{proof}

\begin{prop}\label{p-relative-keel}
Let  $f\colon X \to S$ be a projective morphism of noetherian $\F_p$-schemes. 
Let $L$ be an $f$-nef invertible sheaf  on $X$ and let 
$g\colon \mathbb E_f(L) \hookrightarrow X \xrightarrow{f} S$
be the induced morphism.

Then $\mathbb B_f(L)=\mathbb B_g(L|_{\mathbb E_f(L)})$. In particular, $L$ is $f$-semi-ample if and only if $L|_{\mathbb E_f(L)}$ is $g$-semi-ample. 
\end{prop}

\begin{proof}
Clearly, $\mathbb B_g(L|_{\mathbb E_f(L)})\subset \mathbb B_f(L)$. 
Thus, it is enough to show the opposite inclusion. 
Let $x\in X$ be a 
point such that $x\notin \mathbb B_g(L|_{\mathbb E_f(L)})$. 
Note that if $U\subset S$ is an open subset and if $f'\colon X':=X\times_S U\to U$ is the projection, then $\mathbb E_{f'}(L|_{X'})\subset \mathbb E_f(L)$. 
Thus, we may assume that $S$ is affine.  
By Lemma \ref{l_reduced}, we are reduced to the case where $X$ is reduced.

By (2) of Lemma \ref{l_keel}, we may assume that $L$ is $f$-weakly big. 
Thus, there exist an $f$-ample invertible sheaf $A$ on $X$,  a positive integer $m$ and a nonzero section $s\in H^0(X,L^{\otimes m}\otimes_{\MO_X} A^{-1})$. Let $Z$ be the closed subscheme of $X$ given by the zero set of $s$. Then it follows from 
(1) of Lemma \ref{l_keel} that 
 $\mathbb E_f(L)=\mathbb E_h(L|_Z)$, where $h: Z \hookrightarrow X \xrightarrow{f} S$.  
Since $\Supp Z\subsetneq \Supp X$, 
it follows that
\[
x \notin \mathbb B_g(L|_{\mathbb E_f(L)})=\mathbb B_g(L|_{\mathbb E_h(L|_Z)})
=\mathbb B_h(L|_Z),
\]
where the last equation follows from  noetherian induction. 

We may write $X=X'\cup X''$ where $X'$ (resp. $X''$) is 
the reduced closed subscheme of $X$ 
whose support is equal to the union of all the irreducible components of $X$ that are not contained  (resp. are contained) in $Z$. Thus, $X''\subset Z$ and $D:=X'\cap Z$ is an effective Cartier divisor on $X'$. It follows that $L^{\otimes m}|_{X'}\otimes_{\MO_{X'}} \MO_{X'}(-D)$ is $f'$-ample, where $f'\colon X' \hookrightarrow X \xrightarrow{f} S$ is the induced morphism. 

Since $x\notin \mathbb B_h(L|_Z)$,
  there exist a positive integer $r$ and  $t\in H^0(Z,L^{\otimes mr}|_Z)$ such that $t|_x\neq 0$, where $t|_x$ denotes the pullback of $t$ to $\Spec\,k(x)$ for the residue field $k(x)$ at $x$. 
By Lemma~\ref{l-relative-keel}, there exists a positive integer $e$ and $t'\in H^0(X',L^{\otimes p^emr}|_{X'})$ such that 
\[
t'|_{X'\cap Z}=t^{\otimes p^e}|_{X'\cap Z}.
\]
Since $X''\subset Z$, we have that 
\[
t'|_{X'\cap X''}=t^{\otimes p^e}|_{X'\cap X''}.
\]
By the Mayer--Vietoris type exact sequence 
\[
0 \to \MO_X \to \MO_{X'} \oplus \MO_{X''} \to \MO_{X' \cap X''} \to 0,
\]
we can find a section $u \in H^0(X, L^{\otimes p^emr})$ 
such that $u|_{X'}=t'$ and $u|_{X''}=t^{\otimes p^e}|_{X''}$. 
In particular, $u|_x\neq 0$ and therefore $x\notin \mathbb B_f(L)$. 
Thus, the claim follows. 
\end{proof}

\subsection{Thickening process}

\subsubsection{Partial normalisation}

\begin{dfn}
Let $X$ and $Y$ be reduced noetherian schemes.  
We say that $f\colon Y\to X$ is a {\em partial normalisation} if 
$f$ is a finite birational morphism of schemes. 
In this case,  $Y$ is called a {\em partial normalisation} of $X$. 
\end{dfn}

\begin{dfn}
Let $A$ be a reduced noetherian ring. 
We say that a ring homomorphism $\varphi\colon A \to B$ is 
a {\em partial integral closure} if the induced morphism 
$\Spec\,B \to \Spec\,A$ is a partial normalisation. 
In this case,  $B$ is called a {\em partial integral closure} of $A$. 
\end{dfn}

\begin{rem}
Let $A$ be a reduced noetherian ring whose integral closure 
$A \to A^N$ is finite. 
By definition, a ring homomorphism $\varphi\colon A \to B$ is 
a partial integral closure of $A$ if and only if 
the integral closure $A \to A^N$ factors through $\varphi$. 
If $\varphi\colon A \to B$ is a partial integral closure, then $A$ and $B$ 
admit the same integral closure. 
\end{rem}

\begin{dfn}
Let $A$ be a reduced noetherian ring and let $\varphi\colon A \to B$ be a partial integral closure of $A$. 
We call 
\[
I:=\{a \in A\,|\,aB \subset A\}
\]
the {\em conductor ideal} of $\varphi$. 
Note that $I$ is an ideal of $A$ and also of $B$. 
\end{dfn}

Note that if $A\to B$ is a partial integral closure of a reduced noetherian ring $A$ and $I$ is the conductor ideal, then the sequence 
\[
0 \to A \to B \oplus A/I \to B/I \to 0
\]
is exact, where the third arrow is defined by the difference.

\begin{dfn}
Let $X$ be a reduced noetherian scheme and 
let $f\colon Y \to X$ be a partial normalisation of $X$. 
The closed subschemes $C_X$ and $C_Y$ corresponding to the conductor 
ideals are called {\em conductor subschemes} of $X$ and of $Y$ for $f$, respectively. 
\end{dfn}

\subsubsection{Existence of special thickening subschemes}

\begin{lemma}\label{l-sub-conductor}
Let $A$ be a reduced noetherian ring and let $\varphi\colon A \to B$ 
be a partial integral closure of $A$. 
Let  $I$  be the conductor ideal for $\varphi$. 
Let $J$ be an ideal of $A$ such that $J=JB \cap A$ and $J \subset I$. 

Then the induced sequence 
\[
0 \to A/J \to B/JB \oplus A/I \to B/I \to 0.
\]
is exact, where the third arrow is defined by the difference. 
\end{lemma}

\begin{proof}
The exactness on $A/J$ follows from the assumption $J=JB \cap A$. 
The exactness on $B/I$ is clear. 
The exactness on the middle follows from the fact that the sequence 
$$0 \to A \to B \oplus A/I \to B/I \to 0$$
is exact. 
\end{proof}

\begin{lemma}\label{l-exact-criterion}
Let 
\[
\begin{CD}
A @>\varphi>> B\\
@VV\psi V @VV\psi'V\\
C @>\varphi'>> D
\end{CD}
\]
be a commutative diagram of ring homomorphisms of rings. 
Assume that
\begin{enumerate}
\item $A \to B$ is injective and the induced ring extension is integral. 
\item $A \to C$ is surjective and the above diagram is cocartesian, i.e. 
the induced ring homomorphism $B \otimes_A C \to D$ is bijective. 
\item The sequence 
\[
0 \to A \xrightarrow{(\varphi, \psi)} B \oplus C \xrightarrow{\psi'-\varphi'} D \to 0
\]
is exact.
\end{enumerate}

Then the induced sequence
\[
1 \to \MO_{\Spec\,A}^{\times} \to (\varphi^{\sharp})_*\MO_{\Spec\,B}^{\times} \times 
(\psi^{\sharp})_*\MO_{\Spec\,C}^{\times} \to 
(\psi^{\sharp}\circ \varphi'^{\sharp})_*\MO_{\Spec\,D}^{\times} \to 1.
\]
is exact.
\end{lemma}

\begin{proof}
Fix a prime ideal $\p$ of $A$. 
Let $S:=A \setminus \p$, $A':=S^{-1}A = A_{\p}, B':=S^{-1}B, C':=S^{-1}C$, and $D':=S^{-1}D$. 
Then it is enough to show that the induced sequence 
\[
1 \to A'^{\times} \to B'^{\times} \times C'^{\times} \to D'^{\times} \to 1
\]
is exact. 
After replacing $A$, $B$, $C$ and $D$ by $A'$, $B'$, $C'$, and $D'$ respectively, 
all the assumptions still hold. 
Therefore, 
we may assume that $A$ is a local ring and it suffices to prove that the sequence 
\begin{equation}\label{e1-exact-criterion}
1 \to A^{\times} \to B^{\times} \times C^{\times} \to D^{\times} \to 1
\end{equation}
is exact. 

We first show that $A^{\times}=B^{\times} \cap A$. 
Let $a \in A  \setminus A^{\times}$. 
It suffices to show that $a \not\in B^{\times}$. 
There exists a prime ideal $\p$ of $A$ such that $a \in \p$. 
Since $\Spec\,B \to \Spec\,A$ is surjective by (1), 
there exists a prime ideal $\q$ of $B$ lying over $\p$. 
In particular, we get $a \in \q$, which implies $a \not\in B^{\times}$. 
Thus,  (3) implies that the induced sequence 
\[
1 \to A^{\times} \xrightarrow{(\varphi, \psi)} B^{\times} \times C^{\times} \xrightarrow{\psi'/\varphi'} D^{\times}
\]
is exact. 

In order to prove the exactness of (\ref{e1-exact-criterion}), 
it is enough to show that 
\[
B^{\times}  \to D^{\times}
\]
is surjective. 
Let $I:=\Ker\,\psi$. Then (2) implies that  $D=B/IB$.
Let $d \in D^{\times}$. 
There exist elements $b, b' \in B$ whose images in $D=B/IB$ are equal to $d$ and $d^{-1}$, respectively. 
Thus,  
\[
bb'=1+x
\]
for some $x \in IB \subset \m B$, where $\m$ is the maximal ideal of $A$. 
Since $A \subset B$ is an integral extension by (1), 
\cite[Corollary 5.8]{am69} implies that $\m$ is contained in the Jacobson radical of $B$. 
In particular,  $1+x \in B^{\times}$ and  $b \in B^{\times}$, as desired. 
\end{proof}

\begin{rem}
Note that, using the same notation as in Lemma \ref{l-exact-criterion}, it is easy to check that the condition (3) is equivalent to assuming that $IB=I$, where $I=\Ker\,\psi$. 
\end{rem}

\begin{prop}\label{p-thickening}
Let $f\colon Y \to X$ be a partial normalisation of a reduced noetherian scheme $X$. 
Let $C_X$ and $C_Y$ be the conductor subschemes of $X$ and $Y$, respectively. Let $X_1$ be a closed subscheme of $X$ such that 
$C_X \hookrightarrow X$ factors through $X_1 \hookrightarrow X$. 
Let $Y':=Y \times_X X_1$ and 
let $X'$ be the scheme-theoretic image of $Y'$. 

Then the following hold:
\begin{enumerate}
\item  The closed immersion $C_X \hookrightarrow X_1$ factors through $X'$.
\item $\Supp\, X_1=\Supp\, X'$. 
\item The sequence 
\[
0 \to \MO_{X'} \to \MO_{Y'} \oplus \MO_{C_X} \to \MO_{C_Y} \to 0
\]
is exact, where the third arrow is defined by the difference.
\item 
The sequence 
\[
1 \to \MO^{\times}_{X'} \to \MO^{\times}_{Y'} \times \MO^{\times}_{C_X} \to \MO^{\times}_{C_Y} \to 1
\]
is exact. 
\item Let $L$ be an invertible sheaf on $X$ such that 
\[
L^{\otimes m_1}|_{X'}\simeq \MO_{X'}\qquad\text{and}\qquad  L^{\otimes m_2}|_Y\simeq \MO_Y
\]
for some positive integers $m_1$ and  $m_2$. 
If the restriction map 
\[
H^0(Y, \MO_Y^{\times})_{\Q} \to  H^0(Y', \MO_{Y'}^{\times})_{\Q}
\]
is surjective, then there exists a positive integer $m$ such that  
$L^{\otimes m} \simeq \MO_X.$  
\end{enumerate}
\end{prop}

\begin{proof} 
The assertion (2) follows from the fact that $f$ is proper and surjecitive. 
To prove (1), (3) and (4), we may assume that $X$ and $Y$ are affine: $X=\Spec~A$ and $Y=\Spec~B$. In particular, the induced ring homomorphism $\varphi\colon A\to B$ is a partial integral closure. Let $I$ be the conductor ideal for $\varphi$. 
Let $J_1 \subset A$ and $J' \subset A$ be the ideals of $X_1$ and $X'$, respectively. Since the closed immersion $C_X \hookrightarrow X_1$ 
implies $J_1 \subset I$, we obtain  
$$J'=J_1B \cap A \subset IB \cap A=I \cap A=I.$$
It follows from the definition of $X'$ that $J'=J_1B \cap A$. 
Then 
\cite[Proposition 1.17]{am69} 
implies that $J_1 \subset J'$ and $J'=J'B \cap A$. 
Therefore (1) holds. 
Moreover (3) (resp. (4)) follows from Lemma~\ref{l-sub-conductor} 
(resp. Lemma~\ref{l-exact-criterion}).

We now show (5). Since $C_X$ is contained in $X'$, we have that $L^{\otimes m_1}|_{C_X} \simeq \MO_{C_X}$.

Consider the commutative diagram: 
\[
\begin{CD}
1 @>>> \MO^{\times}_{X} @>>> \MO^{\times}_{Y} \times  \MO^{\times}_{C_X} @>>> \MO^{\times}_{C_Y} @>>> 1\\
@. @VV\alpha V @VV\beta \times {\rm id}V @VV {\rm id} V\\
1 @>>> \MO^{\times}_{X'} @>>> \MO^{\times}_{Y'} \times  \MO^{\times}_{C_X} @>>> \MO^{\times}_{C_Y} @>>> 1,
\end{CD}
\]
where both the horizontal sequences are exact by (4). Thus, we get a commutative diagram 
\[
{\Small
\begin{CD}
H^0(Y,\MO^{\times}_{Y}) \times  H^0(C_X,\MO^{\times}_{C_X}) @>\epsilon >> H^0(C_Y,\MO^{\times}_{C_Y}) @>\delta >> \Pic\,X @>>>\Pic\, Y \times \Pic\, C_X\\
@VV iV @VV {j:=\rm id} V @VV \alpha_1 V @VVV\\
H^0(Y',\MO^{\times}_{Y'}) \times  H^0(C_X,\MO^{\times}_{C_X}) @>\epsilon' >> H^0(C_Y,\MO^{\times}_{C_Y}) @>\delta'>> \Pic\,X' @>>>\Pic\, Y' \times \Pic\, C_X.
\end{CD}
}
\]
where both the horizontal sequences are exact. By a diagram chase, it is easy to check that (5) holds. 
\end{proof}

\subsection{Alteration theorem for quasi-excellent schemes}

The purpose of this subsection is to prove Theorem~\ref{t-alteration}.
Our results  essentially follows from Gabber's alteration theorem for quasi-excellent schemes 
\cite{ILO14}, which in turn  is a generalisation of de Jong's alteration theorem \cite{jong96}. 

We begin by recalling some of the terminology used in \cite{ILO14}. 

\begin{enumerate}
\item[(i)]
A morphism of noetherian schemes $f\colon X \to Y$ is said to be 
{\em generically dominant} 
if the image of any generic point of $X$ by $f$ is a generic point of $Y$ 
\cite[Expos\'e II, D\'efinition 1.1.2]{ILO14}. 
\item[(ii)] Let $S$ be a noetherian scheme. 
We denote by $\alt/S$ the category of reduced $S$-schemes $X$ 
whose structure morphisms $X \to S$ are 
of finite type, generically finite, and generically dominant 
\cite[Expos\'e II, 1.1.9 and D\'efinition 1.2.2]{ILO14}. 
\cite[Expos\'e II, Proposition 1.2.6]{ILO14} implies that the 
category $\alt/S$ admits a fibre product. 
Moreover its proof implies that 
the product of $X$ and $Y$ in $\alt/S$ 
is the reduced closed subscheme 
given by the union of any irreducible component of the scheme-theoretic 
fibre product $X \times_S Y$, which dominates an 
irreducible component of $S$. 
\item[(iii)]
We define the {\em alteration topology}
\cite[Expos\'e II, 2.3.1, 2.3.3]{ILO14}, 
to be the Grothendieck topology on $\alt/S$ 
defined by the pretopology generated by 
\begin{itemize}
\item  \'etale coverings, and 
\item  proper surjective morphisms which are generically finite. 
\end{itemize}
\end{enumerate}

\begin{thm}\label{t-alteration}
Let $X$ be a normal quasi-excellent scheme.
 
Then there exist morphisms of normal quasi-excellent schemes 
$$X_{\nu} \xrightarrow{\varphi_{\nu}} X_{\nu-1} \xrightarrow{\varphi_{\nu-1}} \cdots \xrightarrow{\varphi_{2}} X_1
\xrightarrow{\varphi_{1}} X_0:=X$$
that satisfy the following properties:
\begin{enumerate}
\item 
$X_{\nu}$ is regular. 
\item 
For each $i \in \{1, \cdots, \nu\}$, 
$\varphi_i$ satisfies one of the following: 
\begin{enumerate}
\item $\varphi_i$ is an \'etale surjective morphism.
\item $\varphi_i$ is a morphism which is proper, surjective and generically finite. 
\end{enumerate}
\end{enumerate}
\end{thm}

\begin{proof}
By \cite[Expos\'e II, Th\'eor\`eme 4.3.1]{ILO14} and 
the above definition of the alteration topology, 
there exist morphisms of quasi-excellent reduced schemes 
\[
Y_{\nu} \xrightarrow{\psi_{\nu}} Y_{\nu-1} \xrightarrow{\psi_{\nu-1}} \cdots \xrightarrow{\psi_{2}} Y_1
\xrightarrow{\psi_{1}} Y_0:=X
\]
such that 
\begin{enumerate}
\item[(I)] 
$Y_{\nu}$ is regular. 
\item[(II)] 
For each $i \in \{1, \cdots, \nu\}$, 
one of the following  holds: 
\begin{enumerate}
\item[(A)] $\psi_i$ is an \'etale surjective morphism.
\item[(B)] $\psi_i$ is a morphism which is proper, surjective and generically finite. 
\end{enumerate}
\end{enumerate}
Let $X_i$  be the normalisation of $Y_i$ for each $i$ and let  
\[
X_{\nu} \xrightarrow{\varphi_{\nu}} X_{\nu-1} \xrightarrow{\varphi_{\nu-1}} \cdots \xrightarrow{\varphi_{2}} X_1
\xrightarrow{\varphi_{1}} X_0:=X
\]
be the induced sequence. 
Fix $i \in \{1, \cdots, \nu\}$. 
It is enough to show that (a) or (b) holds. 
Assume (A), i.e. $\psi_i:Y_i \to Y_{i-1}$ is an \'etale surjective morphism. 
Then its base change $\varphi'_i: X_{i-1} \times_{Y_{i-1}} Y_i \to X_{i-1}$ 
is also an \'etale surjective morphism. 
In particular, also $X_{i-1} \times_{Y_{i-1}} Y_i$ is normal. 
Therefore, the induced finite surjective morphism 
$X_{i-1} \times_{Y_{i-1}} Y_i \to Y_i$ coincides with the normalisation. 
Thus, (a) holds. 

If (B) holds, then it is clear that (b) holds. 
\end{proof}

\section{(Theorem~\ref{t-C})$_{n-1}$ implies (Theorem~\ref{t-A})$_n$}

In this section, we prove that (Theorem~\ref{t-C})$_{n-1}$ implies (Theorem~\ref{t-A})$_n$ (cf. Theorem \ref{t_C-to-A}). 
To this end, we first deal with a special case (cf. Proposition \ref{p-normal-nt1}). 
We start with an auxiliary result:

\begin{lemma}\label{l-birat}
Fix a positive integer $n$ and assume $({\rm Theorem}~\ref{t-C})_{n-1}$. 
Let $f\colon X \to S$ be a proper surjective morphism 
of excellent $\F_p$-schemes, 
where $X$ is a normal scheme of dimension $n$. 
Let $L$ be an invertible sheaf on $X$ such that 
$L|_{X_s}$ is semi-ample for all the  points $s \in S$ and 
there exists an open dense subset $S^0$ of $S$ 
such that $L|_{f^{-1}(S^0)}$ is semi-ample over $S^0$ and big over $S^0$.

Then $L$ is $f$-semi-ample.
\end{lemma}

\begin{proof}
We may assume the following properties:
\begin{enumerate}
\item $S$ is an affine scheme. 
\item $X$ and $S$ are integral. 
\item $f_*\MO_X=\MO_S$. In particular $S$ is normal. 
\item $f$ is projective. 
\end{enumerate}
Indeed, we may assume  (1) (resp. (2)) 
by taking an affine open subset (resp. a connected component). By (2) of  Lemma~\ref{l-stein} and  by taking the  Stein factorisation of $f$,  we
may assume (3). 
Finally,  by (4) of Lemma~\ref{l_pullback} and  Chow's lemma,  we may assume (4).

By Proposition~\ref{p-relative-keel}, 
it is enough to show that $L|_{{\mathbb E}_f(L)}$ is relatively semi-ample. By (3) of Lemma~\ref{l_keel}, it follows that $\mathbb E_f(L)$ is a closed subset of $X$.  
Since $S^0$ is a non-empty open subset of $S$ and $L|_{f^{-1}(S^0)}$ is relatively big, 
it follows that $L$ is $f$-weakly big. 
Thus, (2) of Lemma~\ref{l_keel} implies that 
 $\mathbb E_f(L)$ is a proper closed subset of $X$ and, in particular, 
$\dim \mathbb E_f(L)<\dim X$. 
Thus, $({\rm Theorem}~\ref{t-C})_{n-1}$ implies
that 
$L|_{{\mathbb E}_f(L)}$ is relatively semi-ample, as desired.
\end{proof}

\begin{prop}\label{p-normal-nt1}
Fix a positive integer $n$ and assume $({\rm Theorem}~\ref{t-C})_{n-1}$. 
Let $f\colon X \to S$ be a proper morphism 
of excellent $\F_p$-schemes satisfying $f_*\MO_X=\MO_S$, 
where $X$ is a normal scheme of dimension $n$. 
Let $L$ be an invertible sheaf on $X$ such that 
$L|_{X_s}$ is semi-ample for all the points $s \in S$ and 
$L|_{X_\xi}$ is numerically trivial for all the generic points $\xi$ of $S$.

Then $L$ is $f$-semi-ample. 
\end{prop}

\begin{proof}
We may assume that $S$ is affine. 
Replacing $X$ by a purely inseparable model, 
we may assume that the generic fibre of $f$ is geometrically normal.

We want to construct a commutative diagram of morphisms of schemes 
\[
\begin{CD}
X=:X_0 @<\varphi_1 << X_1 @<\varphi_2 << X_2 @<\varphi_3 << \cdots 
@<\varphi_{\nu}<< X_{\nu}\\
@VVf=:f_0 V @VVf_1V @VVf_2V  \cdots @. @VVf_{\nu}V \\
S=:S_0 @<\psi_1 << S_1 @<\psi_2 << S_2@<\psi_3 << \cdots @<\psi_{\nu}<< S_{\nu},
\end{CD}
\]
satisfying the following properties: 
\begin{enumerate}
\item For any $i \in \{1, \cdots, \nu\}$, 
$S_i$ is a normal excellent scheme such that 
$\dim S_i=\dim S$. 
\item For any $i \in \{1, \cdots, \nu\}$, 
$X_i$ is normal excellent schemes such that 
$\dim X_i=\dim X$. 
\item For any $i \in \{1, \cdots, \nu\}$,  $f_i\colon X_i\to S_i$ is a projective surjective morphism such that $(f_i)_*\MO_{X_i}=\MO_{S_i}$. 
\item  For any $i \in \{1, \cdots, \nu\}$ and for any closed point $t \in S_i$, 
we have  $\dim f_i^{-1}(t)=\dim X_i-\dim S_i$. 
\item 
For any $i \in \{1, \cdots, \nu\}$, one of the following holds:
\begin{enumerate}
\item[(a)] $\psi_i$ is an \'etale surjective morphism and 
$X_{i}=X_{i-1} \times_{S_{i-1}} S_i$.
\item[(b)] Both $\varphi_i$ and $\psi_i$ are  proper, surjective and generically finite morphisms,
and $X_{i}$ is the normalisation of the irreducible component of $X_{i-1} \times_{S_{i-1}} S_i$, dominating $S_i$.  
\end{enumerate}
\item $S_{\nu}$ is regular. 
\end{enumerate}
The above diagram can be constructed as follows. 
Below, we denote by $(1)_{i_0}, ..., (5)_{i_0}$ 
the corresponding conditions above in the case $i=i_0$. 

First, $S_1\to S$ is the projective birational morphism so that the projection 
$g_1\colon X'_1 :=X\times_S S_1\to S_1$ is the flattening of $X \to S$, whose existence is guaranteed by \cite[Theorem 5.2.2]{RG71}. 
Let $X_1$  be the normalisation of $X'_1$ and let 
\[
f_1\colon X_1 \to X'_1 \to S_1
\]
be the composite morphism.
Then $(1)_1, \cdots, (5)_1$ hold. 

If $S_1$ is regular, then we are done, otherwise we proceed as follows. The lower horizontal sequence 
\[
S_1 \xleftarrow{\psi_2} S_2 \xleftarrow{\psi_3}  \cdots \xleftarrow{\psi_{\nu}} S_{\nu}
\]
is constructed by applying Theorem~\ref{t-alteration} to $S_1$. 
In particular  $(1)_1,\dots,(1)_\nu$ and $(6)$ hold. 
Moreover, one of the following holds:
 \begin{enumerate}
\item[(a$)'$] $\psi_i$ is an \'etale surjective morphism. 
\item[(b$)'$] $\psi_i$ is a morphism which is proper, surjective and generically finite. 
\end{enumerate}

We now construct $X_i$ inductively as follows. 
Pick $i \in \{1, \cdots, \nu-1\}$ and assume that 
$X_j$, $f_j$ and $\varphi_j$ have already been constructed and 
$(2)_j, \cdots, (5)_j$ hold for any $j\in \{1,\dots,i\}$. 
If $\psi_{i+1}$ satisfies (a$)'$, 
then we define $X_{i+1}:=X_{i} \times_{S_{i}} S_{i+1}$ and 
let $f_{i+1}$ and $\varphi_{i+1}$ be the projections. 
Clearly $(2)_{i+1}, \cdots, (5)_{i+1}$ hold in this case. 
Thus, we may assume that $\psi_{i+1}$ satisfies (b$)'$. 
We provide the construction 
in the case that $X_i, S_i$ and $S_{i+1}$ are integral, 
as we can apply the same argument in the general case 
by taking each connected component separately. 
Since $\psi_{i+1}$ is generically finite, there exists a unique 
irreducible component $X'_{i+1}$ of $X_{i} \times_{S_{i}} S_{i+1}$ 
that dominates $S_{i+1}$, where 
we equip $X'_{i+1}$ with the reduced scheme structure. 
Let $X_{i+1}$ be the normalisation of $X'_{i+1}$. 
Let $f_{i+1}$ and $\varphi_{i+1}$ be the induced morphisms. 
Then  $(2)_{i+1}, (4)_{i+1}$ and $(5)_{i+1}$ hold. 
Further, since $S_{i+1}$ is normal and 
$(f_{i+1})_*\MO_{X_{i+1}}|_{S_{i+1}^0}=\MO_{S_{i+1}}|_{S_{i+1}^0}$ 
for some open dense subset $S_{i+1}^0$ of $S_{i+1}$, also (3)$_{i+1}$ hods. 
This completes the construction of the  commutative diagram above. 

\medskip

For each $i \in \{0, \cdots, \nu\}$,  
the morphism $f_{i}:X_{i} \to S_i$ and the invertible sheaf $L|_{X_i}$ satisfy 
the assumptions in the statement of the proposition. 
We show the claim by descending induction on $i$.

We now show that $L|_{X_{\nu}}$ is $f_{\nu}$-semi-ample. 
To this end, we only treat the case where 
$X_{\nu}$ and $S_{\nu}$ are integral schemes, 
as the general case is reduced to this case by taking connected components. 
By (4)$_{\nu}$ and Lemma~\ref{l-prime-push}, 
it follows that the image of any prime divisor of $X_{\nu}$ 
is either  a prime divisor on $S_{\nu}$ or equal to $S_{\nu}$. 
In particular,  Lemma~\ref{l-flat-descent} implies  
that $L|_{X_{\nu}}$ is $f_{\nu}$-semi-ample. 

Fix $i \in \{0, \cdots, \nu-1\}$ and assume that 
$L|_{X_{i+1}}$ is $f_{i+1}$-semi-ample. 
It is enough to prove that $L|_{X_i}$ is $f_i$-semi-ample. 
If $\psi_{i+1}$ satisfies (a) of (5)$_{i+1}$, then the claim 
follows from (2) of  Lemma~\ref{l-ff-descent}.
Thus, we may assume that $\psi_{i+1}$ satisfies (b) of (5)$_{i+1}$. 
After replacing $X_i, X_{i+1}, S_i$ and $S_{i+1}$ by their connected components, 
we may assume that they are integral schemes. 

We have a commutative diagram: 
\[
\begin{CD}
X_i @<\varphi' << Y@<\varphi''<< X_{i+1}\\
@VVf_iV @VVgV @VVf_{i+1}V\\
S_i @<\psi'<< T@<\psi''<< S_{i+1}\\
\end{CD}
\]
where 
$X_{i+1} \to Y\to  X_i$ is the Stein factorisation of $\varphi_{i+1}$, and 
$Y \to T\to S_i$
is the Stein factorisation of $f_i \circ \varphi'$. 
Note that $S_{i+1} \to S_i$ factors through $T$ because 
$T$ is the Stein factorisation of $X_{i+1} \to S_i$.

Since $L|_{X_{i+1}}$ is $f_{i+1}$-semi-ample, 
it follows from Lemma \ref{l-descend} that there exists a positive integer $m$ and an invertible sheaf $M$ on  $S_{i+1}$ such that 
\[
L^{\otimes m}|_{X_{i+1}} \simeq f_{i+1}^*M.
\]
By Lemma~\ref{l-birat}, $M$ is semi-ample over $T$. Thus,  (1) of Lemma~\ref{l_pullback} implies that
$L|_{X_{i+1}}$ is semi-ample over $T$. 
As $Y$ is normal, (4) of Lemma~\ref{l_pullback}
implies that $L|_Y$ is semi-ample over $T$ and, by (2) of 
Lemma~\ref{l-stein},  it follows that $L|_Y$ is semi-ample over $S_i$.
Since $X_i$ is normal, (4) of Lemma~\ref{l_pullback} implies  that $L|_{X_i}$ is semi-ample over $S_i$. 
This completes the proof.
\end{proof}

\begin{thm}\label{t_C-to-A}
Fix a positive integer $n$. 

Then (Theorem~\ref{t-C})$_{n-1}$ implies  
(Theorem~\ref{t-A})$_{n}$. 
\end{thm}

\begin{proof}
Let $f\colon X \to S$ be a projective surjective morphism 
of excellent $\F_p$-schemes with connected fibres, where $X$ is normal of dimension $n$. 
Let $L$ be an invertible sheaf on $X$ such that 
$L|_{X_s}$ is semi-ample for any point $s \in S$. 
We want to show that $L$ is $f$-semi-ample. 

We may assume the following:
\begin{itemize}
\item $S$ is affine. 
\item $f_*\MO_X=\MO_S$. 
\item $X$ and $S$ are integral normal schemes. 
\end{itemize}
Indeed, we may replace $S$ by an affine open subset. 
By (2) of  Lemma \ref{l-stein}, 
we may replace $f$ by its Stein factorisation. Thus,  we may assume that $f_*\MO_X=\MO_S$ and
in particular, $S$ is normal. 
Replacing $X$ and $S$ by their connected components, 
we may assume that $X$ and $S$ are integral schemes.

We first show the following:
\begin{claim}
There exists a projective birational morphism $\pi\colon Y \to X$ 
and projective morphisms 
\[
g\colon Y \xrightarrow{\varphi} Z \xrightarrow{h} S
\]
of integral normal schemes such that $\varphi_*\MO_Y=\MO_Z$, 
$g=f\circ \pi$ and  $\pi^*L^{\otimes m} = \varphi^*M$, where $m$ is a  positive integer and  $M$ is an invertible sheaf  on $Z$ such that 
$M|_{h^{-1}(S^0)}$ is ample over $S^0$ 
for some open dense subset $S^0$ of $S$. 
\end{claim}

\begin{proof}[Proof of  Claim]
Since $L|_{X_{K(S)}}$ is semi-ample, 
it induces a $K(S)$-morphism 
\[
\psi^1\colon X_{K(S)} \to Z_{K(S)}
\]
to a projective normal $K(S)$-variety $Z_{K(S)}$ with $(\psi^1)_*\MO_{X_{K(S)}}=\MO_{Z_{K(S)}}$. 
Thus, after possibly replacing $L$ by a power of $L$, it follows that  $L|_{X_{K(S)}}$ is the pull-back of an ample invertible sheaf on $Z_{K(S)}$. 

By killing the denominators, 
we can spread out $\psi^1$ over a non-empty open subset $S^0$ of $S$, 
i.e. there exist projective morphisms
\[
f^0\colon X^0=f^{-1}(S^0) \overset{\psi^0}\to Z^0 \overset{h^0}\to S^0
\]
such that $f^0=f|_{f^{-1}(S^0)}$ and
the base change of $\psi^0$ to $K(S)$ is equal to $\psi^1$.  
In particular, $L|_{X^0}$ is the pull-back of an invertible sheaf $M^0$ on $Z^0$ which  is ample over $S^0$. 
Let $Z$ be  a normal projective compactification of $Z^0$ over $S$, so that  we obtain 
\[
X\dashrightarrow Z \xrightarrow{h} S.
\]
Let $Y$ be the normalisation of the resolution of the indeterminacies of $X \dashrightarrow Z$, with induced morphisms  $\pi\colon Y \to X$ and  $\varphi\colon Y\to Z$. Note that $\varphi_*\MO_Y=\MO_Z$.  

Since $\pi^*L|_{Y_z}$ is semi-ample for any $z \in Z$ and $\pi^*L|_{Y_{K(Z)}}$ is numerically trivial, 
 Proposition~\ref{p-normal-nt1} implies that 
$\pi^*L$ is $\varphi$-semi-ample. Thus, Lemma~\ref{l-descend} implies that $\pi^*L^{\otimes m}=\varphi^*M$ where $m$ is a positive integer and $M$ is an invertible sheaf on  $Z$. 
Moreover, after possibly replacing $M^0$ by one of its powers, it follows that  $M|_{h^{-1}(S^0)} \equiv_h M^0$, 
hence $M|_{h^{-1}(S^0)}$ is ample over $S^0$. 
Thus, the claim follows.  
\end{proof}

Since the fibres of $\varphi\colon Y\to Z$ are connected, it follows that also the fibres of 
the restriction morphism $\varphi|_{Y_s}\colon Y_s\to Z_s$ are connected for any $s\in S$. Thus, 
(3) of Lemma~\ref{l_pullback} implies that $M|_{Z_s}$ is semi-ample for any $s \in S$.
 Therefore, Lemma~\ref{l-birat} implies that $M$ is semi-ample over $S$. 
By (1) of Lemma~\ref{l_pullback}, it follows  that   $\pi^*L^{\otimes m} = \varphi^*M$ is semi-ample over $S$. 
Since $X$ is normal, (4) of Lemma~\ref{l_pullback} implies that $L$ is semi-ample over $S$. 
This completes the proof of Theorem~\ref{t_C-to-A}. 
\end{proof}

\section{Numerically trivial case}\label{s-nume-triv}

The main goal of this section is to prove that 
(Theorem \ref{t-A})$_n$ implies (Theorem \ref{t-nume-triv4})$_n$ 
(cf. Theorem \ref{t_A-to-B}). 
In Subsection \ref{ss1-A-to-B}, 
we treat the case where the total space $X$ is normal. 
In Subsection \ref{ss2-A-to-B}, 
we prove that the problem can be reduced to the case where 
the base scheme $S$ is normal. 
In Subsection \ref{ss3-A-to-B}, 
we prove the required statement under the assumption 
that the conductor of the normalisation does not dominate the base scheme.

\subsection{The case where the total space is normal}\label{ss1-A-to-B}

\begin{prop}\label{p-normal-nt}
Fix a positive integer $n$ and assume 
(Theorem~\ref{t-A})$_{n}$.
Let $f\colon X \to S$ be a projective  
morphism of excellent  
$\F_p$-schemes,  
where $X$ is normal of dimension $n$. 
Let $L$ be an $f$-numerically trivial invertible sheaf on $X$ such that 
$L|_{X_s}$ is semi-ample for all the  points $s \in S$. 

Then $L$ is $f$-semi-ample. 
\end{prop}

\begin{proof}
By Lemma~\ref{l-stein}, after possibly taking the Stein factorisation of $f$, we may assume that $f_*\MO_X=\MO_S$. In particular, $S$ is normal. 
Thus, (Theorem~\ref{t-A}$)_{n}$ implies the claim. 
\end{proof}

\subsection{Normalisation of the base}\label{ss2-A-to-B}

We now show that, in order to prove  ({Theorem}~\ref{t-nume-triv4})$_{n}$, we may assume that the base scheme is normal.

\begin{prop}\label{p-reduction-normal}
Fix a positive integer $n$ and assume (Theorem~\ref{t-nume-triv4})$_{n-1}$.
Let 
\[
\begin{CD}
X' @>\alpha >> X\\
@VVf'V @VVf V\\
S' @>\beta >>S
\end{CD}
\]
be a cartesian diagram of excellent $\mathbb F_p$-schemes, 
where $f$ is a projective surjective morphism with connected fibres,  $X$ has dimension $n$ and 
$\beta$ is the composition of the induced morphism $S_{\red} \to S$ and 
the normalisation $S' \to S_{\red}$ of $S_{\red}$.
Let $L$ be an $f$-numerically trivial invertible sheaf on $X$ such that 
$L|_{X_s}$ is semi-ample for all $s \in S$. 
Let $L':=\alpha^*L$. 

Then $L$ is $f$-semi-ample if and only if $L'$ is $f'$-semi-ample. 
\end{prop}

\begin{proof}
By Remark~\ref{r-reduce-hensel}, 
we may assume that $S=\Spec\,R$, where $R$ is a henselian local ring. 
If $L$ is $f$-semi-ample, then (1) of Lemma~\ref{l-stein} and (1) of Lemma~\ref{l_pullback} imply that  $L'$ is  $f'$-semi-ample. 

We now assume that $L'$ is $f'$-semi-ample. 
By Lemma~\ref{l_reduced}, we may assume that $X$ and $S$ are reduced. 
Let $C_S$ and $C_{S'}$ be the conductor subschemes in $S$ and $S'$ 
for $\beta$.
Let $C_X$ and $C_{X'}$  be their inverse images in $X$ and $X'$ respectively. 

\begin{claim} The following hold:
\begin{enumerate}
\item The induced sequence
\[
0 \to \MO_{X}  \to \alpha_*\MO_{X'} \oplus \MO_{C_X}  \to \alpha_*\MO_{C_{X'}} \to 0
\]
is exact. 
\item 
The induced sequence
\[
1 \to \MO^{\times}_{X}  \to \alpha_*\MO^{\times}_{X'} \times  \MO^{\times}_{C_X}  \to \alpha_*\MO^{\times}_{C_{X'}} \to 1
\]
is exact. 
\item 
There exists a positive integer $m_1$ such that $L^{\otimes m_1}|_{X'} \simeq \MO_{X'}$. 
\item 
There exists a positive integer $m_2$ such that $L^{\otimes m_2}|_{C_X} \simeq \MO_{C_X}$. 
\end{enumerate}
\end{claim}

\begin{proof}[Proof of  Claim]
We first show (1). 
By (3) of  Proposition~\ref{p-thickening}, 
we have an exact sequence 
\[
0 \to \MO_{S}  \to \beta_*\MO_{S'} \oplus  \MO_{C_S}  \to \beta_*\MO_{C_{S'}} \to 0.
\]
By Lemma~\ref{l-basic-bc} and by applying $f^*$ to 
the exact sequence above, 
it is enough to show that $\MO_X \to \alpha_*\MO_{X'}$ is injective. 
This follows from the fact that $\alpha\colon X' \to X$ is an affine surjective morphism onto a reduced scheme $X$. 
Thus, (1) holds. Lemma~\ref{l-exact-criterion} implies (2) and Lemma \ref{l-finite/local} implies (3). 
Finally, Lemma \ref{l-finite/local} and 
$({\rm Theorem}~\ref{t-nume-triv4})_{n-1}$ imply (4). 
\end{proof}

By (3) and (4) of 
Claim, after possibly replacing $L$ by $L^{\otimes m_1m_2}$, 
we may assume that $L|_{X'} \simeq \MO_{X'}$ and 
$L|_{C_X} \simeq \MO_{C_X}$.

We have a commutative diagram 
\[
\begin{CD}
1 @>>> \MO^{\times}_{X} @>>> \MO^{\times}_{X'} \times  \MO^{\times}_{C_X} @>>> \MO^{\times}_{C_{X'}} @>>> 1\\
@. @AA\zeta A @AA\xi A @AA \eta A\\
1 @>>> \MO^{\times}_{S} @>>> \MO^{\times}_{S'} \times  \MO^{\times}_{C_S} @>>> \MO^{\times}_{C_{S'}} @>>> 1,
\end{CD}
\]
where, by (2) of  Claim,  both the horizontal sequences are exact. 
Thus, the following diagram is commutative: 
\[
\begin{CD}
H^0(C_{X'},\MO^{\times}_{C_{X'}}) @>\delta_X >> \Pic\,X @>>>\Pic\, X' \times \Pic\, C_X @>>> \Pic\, C_{X'}\\
@AA\eta^0A @AA\zeta^1 A @AA\xi^1A @AA\eta^1 A\\
H^0(C_{S'},\MO^{\times}_{C_{S'}}) @>\delta_S >> \Pic\,S @>>>\Pic\, S' \times \Pic\, C_S @>>> \Pic\, C_{S'}.
\end{CD}
\]
Since $L|_{X'} \simeq \MO_{X'}$ and $L|_{C_X} \simeq \MO_{C_X}$, 
there exists an element $u \in H^0(C_{X'},\MO^{\times}_{C_{X'}})$ 
such that $\delta_X(u) \simeq L$. 
Since $C_{X'} \to C_{S'}$ is a projective morphism with connected fibres, by Lemma~\ref{l_connected-fibres} there is  a positive integer $m$ and an element $v \in H^0(C_{S'},\MO^{\times}_{C_{S'}})$ such that $u^{m}=\eta^0(v)$. 
Therefore, $L^{\otimes m}$ is contained in the image of $\zeta^1$, as desired. 
\end{proof}

\subsection{The vertical case}\label{ss3-A-to-B}

\begin{lemma}\label{l-vertical0}
Fix positive integers $n$ and $m$. 
Assume (Theorem~\ref{t-A})$_{n}$,  
 (Theorem~\ref{t-nume-triv4})$_{n-1}$ and 
(Theorem~\ref{t-nume-triv4})$_{n, m-1}$. 
Let $f\colon X \to S$ be a projective surjective morphism of excellent reduced $\mathbb F_p$-schemes 
with connected fibres, where $X$ has dimension $n$ and 
$S$ is an integral normal scheme of dimension $m$.
Let $L$ be an $f$-numerically trivial invertible sheaf on $X$ 
such that $L|_{X_s}$ is semi-ample for all $s \in S$. 
Assume that
there exists a non-empty open subset $S_1$ of $S$ such that 
the induced morphism $f|_{f^{-1}(S_1)}\colon f^{-1}(S_1) \to S_1$ is a universal homeomorphism. 

Then $L$ is $f$-semi-ample. 
\end{lemma}

\begin{proof}
By Remark~\ref{r-reduce-hensel} and 
the fact that the henselisation of an integrally closed local domain is again an integrally closed local domain, 
we may assume that $S=\Spec\,R$, where $R$ is a henselian local ring. 
We divide the proof into two steps. 

\setcounter{step}{0}

\begin{step}\label{step1-vertical0}
Lemma~\ref{l-vertical0} holds under the assumption that  $X$ is an integral scheme. 
\end{step}

\begin{proof}[Proof of Step~\ref{step1-vertical0}]
In this case, $f\colon X \to S$ is a 
projective surjective morphism of integral excellent schemes. 
By assumption, the induced field extension 
$K(S) \subset K(X)$ is finite and purely inseparable. 
We use the following notation:
\begin{itemize}
\item 
Let $\nu\colon Y \to X$ be the normalisation of $X$. 
Let $C_X$ and $C_Y$ be the conductor subschemes of $X$ and $Y$, respectively. 
Then the composite morphism  
\[
g\colon Y \xrightarrow{\nu} X \xrightarrow{f} S
\]
is a projective surjective morphism of integral normal excellent schemes whose corresponding field extension $K(S) \subset K(Y)$ is finite and purely inseparable. In particular, $g$ has connected fibres. 
\item 
Let $X_1$ be a closed subscheme of $X$ such that the closed immersion 
$C_X \to X$ factors through $X_1$ and 
that 
$\Supp\,X_1$ is equal to $f^{-1}(S')$ where $S':=f(\Supp\,C_X) \cup (S \setminus S_1)$. 
Since $f^{-1}(S_1)\to S_1$ is a universal homeomorphism, it follows that 
$\Supp\,S' \subsetneq S$ and $\Supp\,X_1 \subsetneq X$. 
As $\Supp\,X_1$ is a proper closed subset of a notherian integral scheme $X$, 
it follows that $\dim X_1<\dim X$. 
\item 
Let $Y':=Y \times_X X_1$ and let $X'$ be the scheme-theoretic image of $Y'$. 
By (2) of Proposition~\ref{p-thickening}, it follows that $X'$ and $X_1$ have the same support. In particular, we have that $\dim X'<\dim X$. 
\end{itemize}
By Lemma~\ref{l-finite/local} and (5) of Proposition~\ref{p-thickening}, it is enough to show the following:
\begin{enumerate}
\item[(i)] $L^{\otimes m_1}|_{X'} \simeq \MO_{X'}$  
for some $m_1 \in \Z_{>0}$. 
\item[(ii)] $L^{\otimes m_2}|_{Y} \simeq \MO_Y$  for some $m_2 \in \Z_{>0}$. 
\item[(iii)] 
The restriction map 
\[
H^0(Y,\MO_Y^{\times})_{\Q} \to H^0(Y',\MO_{Y'}^{\times})_{\Q}
\]
 is surjective. 
\end{enumerate}
Thanks to Lemma \ref{l-finite/local}, 
$({\rm Theorem}~\ref{t-nume-triv4})_{n-1}$ implies (i) and, similarly, 
Proposition~\ref{p-normal-nt} implies (ii). 

We now show (iii). 
Note that $\Supp\,Y'=\Supp\,g^{-1}(S')$. 
In particular,  both $g:Y \to S$ and $Y' \to S'$ have connected fibres. Thus, 
by Lemma \ref{l_connected-fibres}, we have the isomorphisms of abelian groups 
\[
H^0(S,\MO_S^{\times})_{\Q} \xrightarrow{\simeq} 
H^0(Y,\MO_Y^{\times})_{\Q}\]
\[
H^0(S',\MO_{S'}^{\times})_{\Q} \xrightarrow{\simeq} 
H^0(Y',\MO_{Y'}^{\times})_{\Q}.
\]
Since $S=\Spec\,R$ where $R$ is a local ring, it follows that  
\[
R^{\times} \to (R/I)^{\times}
\]
is surjective for any ideal $I$ of $R$. 
This implies that 
\[
H^0(S,\MO_S^{\times})_{\Q} \to H^0(S',\MO_{S'}^{\times})_{\Q}
\]
is surjective, hence (iii) holds. 
This completes the proof of Step~\ref{step1-vertical0}. 
\end{proof}

\begin{step}\label{step2-vertical0}
Lemma~\ref{l-vertical0} holds without any additional assumptions. 
\end{step}

\begin{proof}[Proof of Step~\ref{step2-vertical0}]
Let $S_2:=S \setminus S_1$. 
Let $X_1$ be the closure of $f^{-1}(S_1)$ in $X$ and let $X_2:=f^{-1}(S_2)$, where 
we equip $X_1$ and $X_2$ with the reduced scheme structures. 
We denote by $f_1$ the composite morphism: 
\[
f_1\colon X_1 \hookrightarrow X \to S.
\]
The following hold: 
\begin{enumerate}
\item[(I)] $X_1$ and $X_2$ are closed subschemes of $X$. 
\item[(II)] The set-theoretic equality $X=X_1 \cup X_2$ holds. 
\item[(III)] The set-theoretic equality $X_1 \cap X_2=f_1^{-1}(S_2)$ holds. 
\end{enumerate}
By (II) and the fact that $X, X_1$ and $X_2$ are reduced, 
we have the exact sequence 
\[
1 \to \MO_X^{\times} \to \MO_{X_1}^{\times} \times \MO_{X_2}^{\times} \to \MO_{X_1 \cap X_2}^{\times}  \to 1,
\]
which in turn induces the exact sequence 
{\small
\[
H^0(X_1,\MO_{X_1}^{\times})_{\Q} \times H^0(X_2,\MO_{X_2}^{\times})_{\Q} 
\to H^0(X_1\cap X_2,\MO_{X_1 \cap X_2}^{\times})_{\Q} 
\]
\[ 
\to (\Pic\, X)_{\Q} \to (\Pic\, X_1)_{\Q} \times (\Pic\, X_2)_{\Q}.
\]
}
Therefore, it is enough to show the following:
\begin{enumerate}
\item 
$L^{\otimes m_1}|_{X_1} \simeq \MO_{X_1}$ for some $m_1 \in \Z_{>0}$. 
\item 
$L^{\otimes m_2}|_{X_2} \simeq \MO_{X_2}$ for some $m_2 \in \Z_{>0}$. 
\item 
The restriction map 
\[
H^0(X_1,\MO_{X_1}^{\times})_{\Q} \times H^0(X_2,\MO_{X_2}^{\times})_{\Q} \to H^0(X_1\cap X_2,\MO_{X_1 \cap X_2}^{\times})_{\Q}
\]
is surjective. 
\end{enumerate}

By Lemma \ref{l-finite/local}, 
Step \ref{step1-vertical0} implies (1) and 
$({\rm Theorem}~\ref{t-nume-triv4})_{n, m-1}$ implies (2). 

We now show (3). 
Since $X_2=f^{-1}(S_2)$ and $f$ has connected fibres, 
also the induced morphism $X_2 \to S_2$ has connected fibres. 
Thus, Lemma \ref{l_connected-fibres} implies that the induced map 
\begin{equation}\label{e-vertical0-1}
H^0(S_2, \MO_{S_2}^{\times})_{\Q} \to H^0(X_2, \MO_{X_2}^{\times})_{\Q}
\end{equation}
is bijective. 
Since $S$ is normal and $f_1\colon X_1 \to S$ is a proper generically universal homeomorphism of integral schemes, 
it follows that $f_1\colon X_1 \to S$ has connected fibres. Hence, (III) implies that 
 also $X_1 \cap X_2 \to S_2$ has connected fibres. 
Thus, Lemma \ref{l_connected-fibres} implies that 
\begin{equation}\label{e-vertical0-2}
H^0(S_2, \MO_{S_2}^{\times})_{\Q} \to H^0(X_1 \cap X_2, \MO_{X_1 \cap X_2}^{\times})_{\Q}
\end{equation}
is bijective. 
By (\ref{e-vertical0-1}) and (\ref{e-vertical0-2}), we have that 
the map 
\[
H^0(X_2,\MO_{X_2}^{\times})_{\Q} 
\to H^0(X_1\cap X_2,\MO_{X_1 \cap X_2}^{\times})_{\Q} 
\]
is surjective. 
Thus, (3) holds. 
This completes the proof of Step~\ref{step2-vertical0}.
\end{proof}
Step~\ref{step2-vertical0} completes the proof of Lemma~\ref{l-vertical0}. 
\end{proof}

\begin{prop}\label{p-vertical}
Fix positive integers $n$ and $m$. 
Assume (Theorem~\ref{t-A})$_{n}$,  
 (Theorem~\ref{t-nume-triv4})$_{n-1}$ and 
(Theorem~\ref{t-nume-triv4})$_{n, m-1}$. 
Let $f\colon X \to S$ be a projective morphism of excellent reduced 
schemes with connected fibres. 
Let $L$ be an $f$-numerically trivial invertible sheaf on $X$ 
such that $L|_{X_s}$ is semi-ample for all $s \in S$. 
Assume that 
\begin{enumerate}
\item[(a)] 
$\dim X=n$. 
\item[(b)] 
$S$ is an integral scheme. 
\item[(c)] 
The conductor subscheme $C_X$ in $X$ for the normalisation of $X$ satisfies $f(C_X) \subsetneq S$. 
\end{enumerate}

Then $L$ is $f$-semi-ample.
\end{prop}

\begin{proof}
We divide the proof into three steps.

\setcounter{step}{0}
\begin{step}\label{step1-p-vertical}
In order to prove Proposition~\ref{p-vertical}, 
we may assume the following:
\begin{enumerate}
\item 
$S$ is an affine scheme. 
\item
There exists a closed subscheme $\Gamma$ of $X$ such that 
$\Gamma$ is an integral scheme and the induced morphism 
$\Gamma \to S$ is a generically universal homeomorphism. 
\end{enumerate}

\end{step}

\begin{proof}[Proof of Step~\ref{step1-p-vertical}]
Since the problem is local on $S$, we may assume that $S$ is affine.

\begin{claim}
There exists a closed subscheme $T$ of $X$ such that $T$ is 
an integral scheme, 
$T \to S$ is surjective and the induced field extension $K(T) \supset K(S)$ is of finite degree. 
\end{claim}

\begin{proof}[Proof of  Claim]
Take the generic fibre $X \times_S \Spec\,K(S)$, 
which is a scheme of finite type over $K(S)$. 
Since $X \times_S \Spec\,K(S)$ is not empty,  
there exists a closed point $\eta$ of $X \times_S \Spec\,K(S)$. 
It follows from Hilbert's Nullstellensatz that 
$k(\eta) \supset K(S)$ is a finite extension.  
There exists a unique closed subscheme $T$ of $X$ such that 
$T$ is an integral scheme and $T \times_S \Spec\,K(S)$ is equal to $\eta$. 
By construction, $g\colon T \to S$ is dominant. 
Since $g$ is proper, we have that $g$ is surjective. 
It follows from the construction that the field extension 
$K(T) \supset K(S)$ is of finite degree. 
This completes the proof of Claim. 
\end{proof}

Let $L$ be the separable closure of $K(S)$ in $K(T)$. 
By Lemma~\ref{l-local-flattening}, there exists a finite faithfully flat morphism 
\[
S' \to S
\]
where $S'$ is an integral scheme such that $L=K(S')$.  
Take the reduced structure of the base change: 
\[
f'\colon X'=(X\times_S S')_{\red} \to X \times_S S' \to S'.
\]
Clearly the conditions (a) and (b) hold for $X'$ and $S'$. 
Since $S' \to S$ is generically \'etale, 
also the condition (c) holds for $f'$. 
Since $S' \to S$ is faithfully flat, 
we can replace $f$ by $f'$ (Lemma~\ref{l-ff-descent}). 
By  construction, we can find the required closed subscheme 
$\Gamma$ of $X'$ as an irreducible component of $S' \times_S T$. 
This completes the proof of Step~\ref{step1-p-vertical}. 
\end{proof}

\begin{step}\label{step2-p-vertical}
In order to prove Proposition~\ref{p-vertical}, 
we may assume 
the condition (2) in  Step~\ref{step1-p-vertical} 
and the following conditions (3) and (4). 
\begin{enumerate}
\item[(3)]
$S$ is normal. 
\item[(4)]
$S=\Spec\,R$, where $R$ is a henselian local ring.  
\end{enumerate}
\end{step}

\begin{proof}[Proof of Step~\ref{step2-p-vertical}]
By Step~\ref{step1-p-vertical}, 
we may assume that $f\colon X \to S$ satisfies (1) and (2). 
Let $S' \to S$ be the normalisation of $S$ and 
consider the reduced structure of the base change 
\[
f'\colon X'=(X\times_S S')_{\red} \to X \times_S S' \to S'.
\]
Clearly, (a), (b), (c), (1), (2) and (3) hold for $f'\colon X' \to S'$. 
By Proposition~\ref{p-reduction-normal}, 
we may replace $f$ by $f'$. 
Thus, we may assume that (1), (2) and (3) hold. 
By Remark~\ref{r-reduce-hensel}, 
we are done. 
Note that the henselisation does not break the condition (b) in our case. 
Indeed, if $R$ is a normal excellent local ring, then 
so is $R^h$, hence in particular $R^h$ is an integral domain. 
\end{proof}

\begin{step}\label{step3-p-vertical}
Proposition \ref{p-vertical} holds without any additional assumptions.
\end{step}

\begin{proof}[Proof of Step \ref{step3-p-vertical}]
By Step~\ref{step2-p-vertical}, we may assume that 
(2)--(4) hold. 
Let $\nu\colon Y \to X$ be the normalisation of $X$. 
Let $g\colon Y \to T$ be the Stein factorisation of $Y \to S$. 
We can find a closed subscheme $S_1$ of $S$ such that 
\begin{itemize}
\item $\Supp\, S_1 \subsetneq \Supp\, S$,  
\item $f(C_X) \subset \Supp\,S_1$, and 
\item $\Gamma \setminus f^{-1}(S_1) \to S \setminus S_1$ is 
a universal homeomorphism. 
\end{itemize}
We take a closed subscheme $X_1$ of $X$ 
such that $\Supp \, X_1=\Gamma \cup f^{-1}(S_1)$. 
Let $Y':=Y \times_X X_1$ and let $X'$ be the scheme-theoretic image of $Y'$. 
By  (2) of Proposition~\ref{p-thickening}, it follows that $X'$ and $X_1$ have the same support. 
It follows that $X' \to S$ and $Y' \to T$ have connected fibres. 

By Lemma~\ref{l-finite/local} and (5) of Proposition~\ref{p-thickening}, it is enough  to show the following: 
\begin{enumerate}
\item[(i)] $L|_{X'}$ is semi-ample over $S$. 
\item[(ii)] $L|_{Y}$ is semi-ample over $S$. 
\item[(iii)] The restriction map 
\[
H^0(Y, \MO_Y^{\times})_{\Q} \to  H^0(Y', \MO_{Y'}^{\times})_{\Q}
\]
is bijective. 
\end{enumerate}
Thanks to (Theorem~\ref{t-nume-triv4})$_{n-1}$ and (Theorem~\ref{t-nume-triv4})$_{n, m-1}$, 
we may apply Lemma~\ref{l-vertical0}, hence (i) holds. 
Proposition~\ref{p-normal-nt} implies (ii). 
Since both the morphisms $Y \to T$ and $Y' \to T$ have connected fibres, Lemma~ \ref{l_connected-fibres} 
implies (iii). 
This completes the proof of Step \ref{step3-p-vertical}. 
\end{proof} 
Step \ref{step3-p-vertical} completes the proof of Proposition~\ref{p-vertical}. 
\end{proof}

\subsection{ ({Theorem~\ref{t-A}})$_{n}$  implies  
(Theoerem~\ref{t-nume-triv4})$_n$}
\label{s_num-triv}

\setcounter{step}{0}

\begin{thm}\label{t_A-to-B}
Fix a positive integer $n$. 

Then   (Theorem~\ref{t-A})$_{n}$ implies 
(Theorem~\ref{t-nume-triv4})$_n$.
\end{thm}

\begin{proof}
We first introduce some notation. 

Let $f\colon X \to S$ be as in the statement of Theorem~\ref{t-nume-triv4}. 
Let $m=\dim S$ and let $S_1, \cdots, S_t$ be the $m$-dimensional irreducible components of $S$ 
equipped with the reduced scheme structures. 
For any $k \in \{1, \cdots, t\}$, let $\xi_k$ be the generic point of $S_k$ 
and let $\overline{\xi}_k$ be the geometric point obtained by taking 
its algebraic closure. 
Let
\[
\delta(f):=\max_{1 \leq k \leq t} \dim X_{\xi_k}.
\]

Let $\nu\colon X^N \to X$ be the normalisation of $X$ and 
let $C_X$ be the conductor subscheme of $X$ for $\nu$. 
For any $k \in \{1, \cdots, t\}$, 
let $\eta_k(f)$  be the number of the connected components of 
the fibre $C_{X, \overline{\xi}_k}$ over $\overline{\xi}_k$ of the induced morphism
\[
C_X \hookrightarrow X \to S.
\]
Let
\[
\eta(f):=\max_{1 \leq k \leq t} \eta_k(f).
\]
We consider the   set-theoretic decomposition 
\[
C_X=C_X^h \cup C_X^v
\]
so that $C_X^h$ and $C_X^v$ are closed subsets of $X$ 
which admit decompositions into irreducible components 
\[
C_X^h=\bigcup_{i=1}^r C_X^{h, i}, \quad C^v_X =\bigcup_{j=1}^s C_X^{v, j}
\]
as closed subsets of $X$, 
where each $C_X^{h, i}$ dominates $S_k$ for some $k \in \{1, \cdots, t\}$ 
and each $C_X^{v, j}$ does not dominate any of $S_1, \cdots, S_t$. 
We equip $C_X^{h, i}$ and $C_X^{v, j}$ with the reduced scheme structures. 
In particular, each of $C_X^{h, i}$ and $C_X^{v, j}$ is an integral scheme. 
 
Let 
\[
Q(f):=(\dim X, \dim S,  \delta(f), \eta(f)).
\]
We proceed by induction on all the quadruples of non-negative integers $(n,m,\delta,\eta)$ with respect to  the lexicographic order 
(e.g. $(1, 0, 0, 0)>(0, 1, 0, 0)$).

\begin{step}\label{step1-p-nume-triv}
Let $f\colon X \to S$ be as in the statement of Theorem~\ref{t-nume-triv4}. 
Let  $f'\colon X \to S'$ be the Stein factorisation of $f\colon X \to S$. 
Then the following hold:
\begin{itemize}
\item $S'$ is reduced.
\item $\dim S=\dim S'$. 
\item $\delta(f) = \delta(f')$. 
\item $\eta(f) \geq \eta(f')$. 
\end{itemize}
In particular, $Q(f)\ge Q(f')$. 
\end{step}

\begin{proof}[Proof of Step~\ref{step1-p-nume-triv}]
Let $\beta:S' \to S$ be the induced morphism. 
Let $S'_1, \cdots, S'_{t'}$ be the $m$-dimensional irreducible componenets of $S'$ 
and let $\xi'_{\ell}$ be the generic point of $S'_{\ell}$ for $\ell \in \{1, \cdots, t'\}$. 
Since $X$ is reduced, so is $S'$. 
Note that for any open affine subset $\Spec\,R$ of $S$, if we denote by $\Spec\,R'$ its inverse image to $S'$, 
then  $R \to R'$ is a finite injective ring homomorphism. 
Thus, \cite[Theorem 5.11]{am69} implies the following hold: 
\begin{itemize}
\item $\dim S=\dim S'$. 
\item For any $\ell \in \{1, \cdots, t'\}$, there exists $k \in \{1, \cdots, t\}$ 
such that $\beta(\xi'_{\ell})=\xi_k$. 
\item For any $k \in \{1, \cdots, t\}$, there exists a non-empty 
subset $L$ of $\{1, \cdots, t'\}$ such that $\beta^{-1}(\{\xi_k\})=\bigcup_{\ell \in L}\{\xi'_{\ell}\}$. 
\end{itemize}
Thus, it follows that 
$\delta(f) = \delta(f')$ and $\eta(f) \geq \eta(f')$. 
This completes the proof of Step~\ref{step1-p-nume-triv}.
\end{proof}

\begin{step}\label{step2-p-nume-triv}
Let $f\colon X \to S$ and let $L$ be as in the statement of Theorem~\ref{t-nume-triv4}. 
Let 
\[
\beta\colon S'' \to S
\]
be a morphism satisfying one of the following properties:
\begin{itemize}
\item $\beta$ is the normalisation of $S$. 
\item $S$ and $S''$ are integral schemes and 
$\beta$ is a finite flat generically \'etale morphism. 
\end{itemize}
Consider the reduced structure of the base change of $f$ over $S''$: 
\[
f''\colon X''=(X \times_S S'')_{\red} \to X \times_S S'' \to S''.
\]
Then the following hold:
\begin{itemize}
\item $\dim X=\dim X''$. 
\item $\dim S=\dim S''$. 
\item $\delta(f) = \delta(f'')$. 
\item $\eta(f) = \eta(f'')$. 
\item If $L|_{X''}$ is $f''$-semi-ample, then $L$ is $f$-semi-ample. 
\end{itemize}
In particular, $Q(f)=Q(f'')$. 
\end{step}

\begin{proof}[Proof of Step~\ref{step2-p-nume-triv}] 
Since $\beta$ is a finite surjective morphism, so is $X'' \to X$. 
Thus, we have that $\dim S=\dim S''$ and $\dim X=\dim X''$. 
It is easy to check that 
$\delta(f) = \delta(f'')$ and $\eta(f) = \eta(f'')$. 
By Lemma~\ref{l-ff-descent} and Proposition~\ref{p-reduction-normal}, we have that if $L|_{X''}$ is $f''$-semi-ample, then $L$ is $f$-semi-ample.
This completes the proof of Step~\ref{step2-p-nume-triv}.
\end{proof}

\begin{step}\label{step3-p-nume-triv} Fix positive integers $n$ and $m$. 
Assume (Theorem~\ref{t-nume-triv4})$_{n-1}$ and (Theorem~\ref{t-nume-triv4})$_{n, m-1}$. 

Then 
(Theorem~\ref{t-nume-triv4})$_{n, m}$ 
 holds for any morphism $f\colon X\to S$ such that $\delta(f)=0$ or $\eta(f)=0$. 
\end{step}

\begin{proof}[Proof of Step~\ref{step3-p-nume-triv}]
Let $f\colon X \to S$ and $L$ be as in 
$({\rm Theorem}~\ref{t-nume-triv4})_{n, m}$ and such that $\delta(f)=0$ or $\eta(f)=0$. 
By Step~\ref{step1-p-nume-triv} and Step~\ref{step2-p-nume-triv}, 
we may assume that $f$ has connected fibres and $S$ is normal. 
Since the problem is local on $S$, we may assume that $S$ is an integral normal scheme. 
Since $\delta(f)=0$ or $\eta(f)=0$, we have that 
$f(C_X) \subsetneq S$ where $C_X$ denotes the conductor of the normalisation of $X$. 
Thus, 
Proposition~\ref{p-vertical} 
implies that  $L$ is $f$-semi-ample. 
\end{proof}

\begin{step}\label{step5-p-nume-triv} Fix positive integers $n$, $m$, $\delta$ and $\eta$. Assume 
 that Theorem~\ref{t-nume-triv4} holds for all the morphisms $f\colon X\to S$ such that $Q(f)<(n,m,\delta,\eta)$.

Then Theorem~\ref{t-nume-triv4} 
 holds for any morphism $f\colon X\to S$ such that  $Q(f)=(n,m,\delta,\eta)$
 and satisfying the following properties: 
\begin{enumerate}
\item[(a)] $f\colon X \to S$ has connected fibres. 
\item[(b)] $S=\Spec\,R$, where $R$ is an integral normal local henselian ring. 
\item[(c)] The induced morphism $f^{h, 1}\colon C_X^{h, 1} \to S$ has connected fibres.  
\end{enumerate}
\end{step}

\begin{proof}[Proof of Step~\ref{step5-p-nume-triv}]
Let $\nu\colon X^N \to X$ be the normalisation of $X$. 
By \cite[Theorem 7.1]{ferrand03}, 
we can find morphisms 
\[
\nu\colon X^N \to Y \xrightarrow{\pi_1} X
\]
such that 
\begin{enumerate}
\item[(i)] $Y$ is a reduced scheme. 
\item[(ii)] Both $X^N \to Y$ and $Y\to X$ are finite birational morphisms. 
\item[(iii)] The conductor $D_X$ of $X$ for $\pi_1\colon Y\to X$ is set-theoretically equal to $C_X^{h, 1}$. 
\item[(iv)] If $f_{Y}\colon Y \to X \to S$ is the induced morphism, then $\eta(f)>\eta(f_{Y})$. 
\item[(v)] Any irreducible component of the conductor $D_Y$ of $Y$ for $\pi_1\colon Y \to X$ dominates $S$.  
\end{enumerate}

Indeed, such $Y$ can be constructed as follows. 
If $C'_X$ denotes the scheme-theoretic image of 
the induced immersion $C_X \cap (X \setminus C_X^{h, 1}) \to X$, 
then we define $Z_1$ as the pushout of the diagram $X^N \hookleftarrow \nu^{-1}(C'_X) \to C'_X$, whose existence is guaranteed by \cite[Theorem 7.1]{ferrand03}. 
Let $Z:=(Z_1)_{\red}$. 
Then $Z$ satisfies the corresponding properties (i)$_Z$--(iv)$_Z$ to (i)--(iv). 
We denote by $E_X$ and $E_Z$  the conductors of $X$ and $Z$ respectively 
for the induced finite birational morphism $\mu \colon Z \to X$. 
Let 
\[
\Supp\,E_Z=(E_1 \cup \cdots \cup E_a) \cup (F_1 \cup \cdots \cup F_b)
\]
be the irreducible decomposition such that 
all of $E_1, \cdots, E_a$ dominate $S$ and 
none of $F_1, \cdots, F_b$ dominates $S$. 
Let $C''_X$ be the reduced closed subscheme of $X$ 
that is set-theoretically equal to $\mu(F_1 \cup \cdots \cup F_b)$. 
We define $Y$ as the pushout of the diagram $Z \hookleftarrow \mu^{-1}(C''_X) \to C''_X$, whose existence is guaranteed again by \cite[Theorem 7.1]{ferrand03}. 
Since $Z$ and $C''_X$ are reduced, so is $Y$. Hence (i) holds. 
The properties (ii), (iii) and (v) follow directly from the construction. 
The remaining one (iv) holds by (iv)$_Z$ and the fact that 
the induced morphism $Z_{K(S)} \to Y_{K(S)}$ of the generic fibres is an isomorphism.

\medskip

Let $\eta =\Spec~K(S)$ be the generic point of $S$ and let $X_\eta=X\times_S \Spec\, K(S)$ be the generic fibre of $f$. Similarly, we denote 
\[
Y_\eta=Y\times_S \Spec\, K(S), \quad D_{X_\eta}=D_{X}\times_S \Spec \,K(S)\quad \text{and} \quad D_{Y_\eta}=D_Y \times_S \Spec \, K(S).
\]
Note that $D_{X_\eta}$ and $D_{Y_\eta}$ are the conductor of the morphism $Y_\eta\to X_\eta$ in $X_\eta$ and $Y_\eta$ respectively. 
We have the commutative diagram:
\[
{\Small
\begin{CD}
H^0(Y,\MO^{\times}_{Y}) \times  H^0(D_X,\MO^{\times}_{D_X}) @>\varphi >> H^0(D_Y,\MO^{\times}_{D_Y}) @>\psi >> \Pic\,X @>>>\Pic\, Y \times \Pic\, D_X\\
@VV iV @VV {j} V @VVV @VVV\\
H^0(Y_\eta,\MO^{\times}_{Y_{\eta}}) \times  H^0(D_{X_\eta},\MO^{\times}_{D_{X_\eta}}) @>\varphi' >> H^0(D_{Y_\eta},\MO^{\times}_{D_{Y_\eta}}) @>\psi'>> \Pic\,X_{\eta} @>>>\Pic\, Y_{\eta} \times \Pic\, D_{X_\eta}.
\end{CD}
}
\]

\begin{claim} 
The following hold: 
\begin{enumerate}
\item There exists a positive integer $r$ 
such that 
\[
L^{\otimes r}|_Y \simeq \MO_Y, \quad L^{\otimes r}|_{D_X} \simeq \MO_{D_X}
\quad{\rm and}\quad L^{\otimes r}|_{X_{\eta}} \simeq \MO_{X_{\eta}}.
\]
\item ${\rm Im}(\varphi'_{\Q}) \cap {\rm Im}(j_{\Q}) \subset {\rm Im}(\varphi'_{\Q} \circ i_{\Q}).$ 
\item $j_{\Q}\colon H^0(D_Y,\MO^{\times}_{D_Y})_{\Q} \to H^0(D_{Y_\eta},\MO^{\times}_{D_{Y_\eta}})_{\Q}$ is injective. 
\end{enumerate}
\end{claim}

\begin{proof}[Proof of Claim]
We first show (1). 
Since we are assuming that 
Theorem~\ref{t-nume-triv4} holds for all the morphisms $f$ such that $Q(f)<(n,m,\delta,\eta)$, we have that $L|_Y$ and $L|_{D_X}$ are semi-ample over $S$. 
Thus, Lemma~\ref{l-finite/local} implies that there exist $r_1 \in \Z_{>0}$ such that 
$L^{\otimes r_1}|_{Y}\simeq \MO_Y$ and $L^{\otimes r_1}|_{D_X}\simeq \MO_{D_X}$. 
By assumption, $L|_{X_{\eta}}$ is semi-ample. 
Hence, again by Lemma~\ref{l-finite/local}, we may find $r_2 \in \Z_{>0}$ 
such that $L^{\otimes r_2}|_{X_{\eta}} \simeq \MO_{X_{\eta}}$. 
Let $r:=\max\{r_1,r_2\}$. Then (1) holds.

We now show (2). 
Let 
\[
R_Y:=H^0(Y,\MO_{Y}), \quad R_{D_X}:=H^0(D_X, \MO_{D_X}), 
\quad R_{D_Y}:=H^0(D_Y, \MO_{D_Y}).
\]
Since these rings define the Stein factorisations of $Y \to S$, $D_X \to S$, and $D_Y \to S$ respectively, we obtain injective ring homomorphisms
\begin{equation}\label{C.1}
\Gamma(S, \MO_S)=R \to R_Y \to R_{D_Y}. 
\end{equation}
For $U:=R \setminus \{0\}$, 
the left square in the  diagram above induces the commutative diagram
\[
\begin{CD}
(R_Y^{\times})_{\Q} \times (R_{D_X}^{\times})_{\Q} @>\varphi_{\Q}>>(R_{D_Y}^{\times})_{\Q}\\
@VVi_{\Q}V @VVj_{\Q}V\\
(U^{-1}R_Y)^{\times}_{\Q} \times (U^{-1}R_{D_X})^{\times}_{\Q}  @>\varphi'_{\Q}>>
(U^{-1}R_{D_Y})^{\times}_{\Q}. \\
\end{CD}
\]
Since $D_X \to S$ has connceted fibres, 
 Lemma \ref{l_connected-fibres} implies that 
\begin{equation}\label{C.2}
(U^{-1}R_{D_X}^{\times})_{\Q}=(K(S)^{\times})_{\Q}.
\end{equation}
Pick $\alpha \in {\rm Im}(\varphi'_{\Q}) \cap {\rm Im}(j_{\Q})$. 
We want to show that  $\alpha \in {\rm Im}(\varphi'_{\Q} \circ i_{\Q})$. 
It follows from (\ref{C.1}) and (\ref{C.2}) that, 
possibly after replacing $\alpha$ by $\alpha^s$ for some $s \in \Z_{>0}$, 
there exist 
$\beta \in (U^{-1}R_Y)^{\times}$ and $\gamma \in R_{D_Y}^{\times}$ such that 
\[
\alpha=\varphi'(\beta, 1)=j(\gamma). 
\]
In particular, we have that $\beta, \beta^{-1} \in K(R_Y) \cap R_{D_Y}^\times$. 
Since $R_Y$ is an integrally closed integral domain and $R_Y \to R_{D_Y}$ 
is a finite injective ring homomorphism, 
\cite[Theorem 9.1]{matsumura89} implies 
that $\beta, \beta^{-1} \in R_Y$. 
In particular, we get $\beta \in R_Y^{\times}$. 
It follows that 
\[
\alpha=\varphi'(\beta, 1)=(\varphi' \circ i)(\beta, 1) \in 
{\rm Im}(\varphi'_{\Q} \circ i_{\Q}).
\]
Thus, (2) holds.

Finally, we show (3).  
Let $(D_Y)_{\red}^N$ be the normalisation of $(D_Y)_{\red}$. 
We have a commutative diagram: 
\[
\begin{CD}
H^0((D_Y)^N_{\red},\MO^{\times}_{(D_Y)^N_{\red}})_{\Q} @>(j^N_{\red})_{\Q}>> H^0((D_{Y_{\eta}})^N_{\red},\MO^{\times}_{(D_{Y_{\eta}})^N_{\red}})_{\Q}\\
@AA\nu A @AA\nu_{\eta}A\\
H^0((D_Y)_{\red},\MO^{\times}_{(D_Y)_{\red}})_{\Q} @>(j_{\red})_{\Q}>> H^0((D_{Y_{\eta}})_{\red},\MO^{\times}_{(D_{Y_{\eta}})_{\red}})_{\Q}\\
@AA\rho A @AA\rho_{\eta} A\\
H^0(D_Y,\MO^{\times}_{D_Y})_{\Q} @>j_{\Q}>> H^0(D_{Y_\eta},\MO^{\times}_{D_{Y_\eta}})_{\Q}.
\end{CD}
\]
Clearly, $\nu$ is injective. 
Since any irreducible component of $D_Y$ dominates $S$, 
it follows that $(j^N_{\red})_{\Q}$ is injective. 
 Lemma \ref{l_connected-fibres} implies that $\rho$ is bijective. 
Therefore, the composite map $(j^N_{\red})_{\Q} \circ \nu \circ \rho$ 
is injective and, in particular,  also $j_{\Q}$ is injective. 
Thus, (3) holds. 
\end{proof}

By (1) of Claim, after possibly replacing $L$ by one of its powers,  
we may assume that 
$L|_Y \simeq \MO_Y, L|_{D_X} \simeq \MO_{D_X}$ 
and $L|_{X_{\eta}} \simeq \MO_{X_{\eta}}.$ 
Thus there exists an element $a\in H^0(D_Y,\MO^{\times}_{D_Y})$ such that $\psi(a) = L$. 
Since $L|_{X_{\eta}} \simeq \MO_{X_{\eta}},$ 
it follows that $j(a) \in {\rm Im}(\varphi')$. 
Hence, after possibly replacing $L$ again, 
by (2) of Claim, we may assume that there exists an element 
$b \in  H^0(Y,\MO^{\times}_{Y}) \times  H^0(D_X,\MO^{\times}_{D_X})$ 
such that $(\varphi' \circ i)(b)=j(a)$ and, in particular, $j(a\cdot \varphi(b)^{-1})=1$. 
By (3) of Claim, it follows that $a^s=\varphi(b^s)$ for some $s \in \Z_{>0}$. 
This implies that $L^{\otimes s} = \psi(a^s)= (\psi \circ \varphi)(b^s) = \MO_X$. 
This completes the proof of Step~\ref{step5-p-nume-triv}. 
\end{proof}

\begin{step}\label{step6-p-nume-triv}
Fix positive integers $n$, $m$, $\delta$ and $\eta$. Assume 
 that Theorem~\ref{t-nume-triv4} holds for all the morphisms $f\colon X\to S$ such that $Q(f)<(n,m,\delta,\eta)$.

Then Theorem~\ref{t-nume-triv4} 
 holds for all the morphisms $f\colon X\to S$ such that  $Q(f)=(n,m,\delta,\eta)$.
\end{step}

\begin{proof}[Proof of Step~\ref{step6-p-nume-triv}]
We shall reduce the problem to the case where (a), (b) and (c) in 
Step~\ref{step5-p-nume-triv} hold. 
By Lemma~\ref{l-stein},  Step~\ref{step1-p-nume-triv}, Step~\ref{step2-p-nume-triv} 
and the fact that the problem is local on $S$, 
we may assume the following:
\begin{enumerate}
\item[(a)] $f\colon X \to S$ has connected fibres. 
\item[(b)'] $S$ is an integral normal affine scheme. 
\end{enumerate}
By Lemma~\ref{l-local-flattening}, 
we can find a finite faithfully flat separable morphism $S_1 \to S$ 
from an integral affine scheme $S_1$ such that 
for the commutative diagram 
\[
\begin{CD}
X_1:=(X \times_S S_1)_{\red} @>\alpha >> X\\
@VVf_1 V @VVfV\\
S_1 @>\beta >> S,
\end{CD}
\]
there exists an irreducible component $D$ of $\alpha^{-1}(C_X^{h, 1})$ 
equipped with the reduced scheme structure 
such that $D \to S_1$ has connected fibres. 
Since $S_1 \to S$ is separable i.e. generically \'etale, 
we have that $\alpha^{-1}(C_X^h)$ and the horizontal part $C_{X_1}^h$ 
 coincide over the open subset of $S_1$ where $\beta$ is \'etale. 
In particular, (a) and (c) holds for $f_1$. 

Let $S_2$ be the normalisation of $S_1$ and set 
\[
f_2\colon X_2=(X_1 \times_{S_1} S_2)_{\red} \to S_2.
\]
By Step~\ref{step2-p-nume-triv}, 
we may replace $f$ by $f_2$. 
In particular, $f$ satisfies (a), (b)' and (c). 
Finally replacing $X \to S$ by the base change of the henselisation 
of a stalk of $S$, 
we may assume (b). 
This completes the proof of Step~\ref{step6-p-nume-triv}. 
\end{proof}

By quadruple induction on $n$, $m$, $\delta$ and $\eta$, it follows that
Step~\ref{step3-p-nume-triv} and 
Step~\ref{step6-p-nume-triv} complete the proof of Theorem \ref{t_A-to-B}. 
\end{proof}

\subsection{Generalisation to algebraic spaces}

\begin{thm}\label{t-nume-triv-as}
Fix a positive integer $n$ and assume $({\rm Theorem}~\ref{t-nume-triv4})_{n}$.
Let $f\colon X \to S$ be a projective surjective morphism 
of excellent algebraic spaces over $\F_p$ 
with connected fibres, where $X$ is of dimension $n$. 
Let $L$ be an invertible sheaf on $X$ such that 
$L$ is $f$-numerically trivial and $L|_{X_s}$ is semi-ample 
for all the  points $s$ of $S$. 

Then there exists a positive integer $m$ 
and an invertible sheaf $M$ on $S$ such that 
\[
L^{\otimes m} \simeq f^*M.
\]
\end{thm}

\begin{proof}
Let $f\colon X\xrightarrow{g} T \xrightarrow{\eta} S$ 
be the Stein factorisation of $f$. 
Since the fibres of $f$ are connected, 
$\eta$ is a finite universal homeomorphism.
By \cite[Proposition 6.6]{kollar97}, there exists a positive integer $e$ such that 
the $e$-th iterated Frobenius morphism $F^e\colon T\to T$ factors through $\eta$. 
Therefore, replacing $f$ by $g$, we are reduced the case where $f_*\MO_X=\MO_S$. 

Let $\beta\colon S' \to S$ be an \'etale surjective morphism from a scheme $S'$. 
Let $X':=X \times_S S'$, so that the following diagram is cartesian:
\[\begin{CD}
X' @>\alpha >> X \\
@VVf' V @VVf V\\
S' @>\beta >> S.
\end{CD}
\]
Since  the induced morphism $f'\colon X' \to S'$ is projective, 
it follows that $X'$ is a scheme (cf. \cite[Ch. II, Definition 7.6]{knutson71}). 
Therefore, (Theorem~\ref{t-nume-triv4})$_n$ implies that $L|_{X'}$ is semi-ample over $S'$. 
After possibly replacing $L$ by one of its powers, 
we have that   
\[
\alpha^*L \simeq f'^*N
\] 
for some invertible sheaf $N$ on $S'$.

We now show that $f_*(L)$ is an invertible sheaf on $S$. 
By \cite[Ch. II, Definition 4.1]{knutson71}, it is enough to show that 
$\beta^*(f_*L))$ is an invertible sheaf. 
By the flat base change theorem, we have that 
\[
f'_*(\alpha^*L) \simeq \beta^*(f_*L)).
\]
Since $f'_*\MO_{X'}=\MO_{S'}$ and $\alpha^*L \simeq f'^*N$, 
we have that 
\[
f'_*(\alpha^*L) \simeq f'_*(f'^*N) \simeq N.
\]
Hence, $\beta^*(f_*L))$ is an invertible sheaf, as desired. 

We have that the induced homomorphism 
\[
\theta\colon f^*f_*L \to L
\]
is surjective, since so is its pull-back by $\alpha$. 
Since both $f^*f_*L$ and $L$ are invertible, 
it follows from \cite[Theorem 2.4]{matsumura89} 
that $\theta$ is an isomorphism. 
\end{proof}

\section{
(Theorem~\ref{t-A}$)_n$ and $({\rm Theorem~\ref{t-nume-triv4}})_n$ imply (Theorem~\ref{t-C}$)_n$ }\label{s-AB-to-C}

\begin{dfn}[Definition 9.1, 9.2 and 9.4 of \cite{kollar13}]\label{d-quotients}
Let $S$ be a scheme and let $X$ be an algebraic space over $S$. 
\begin{enumerate}
\item 
A {\em relation} on $X$ over  $S$ is a closed immersion $\sigma\colon R \to X \times_S X$ over $S$, where 
$R$ is an algebraic space  over $S$.
\item 
A relation $\sigma\colon R \to X \times_S X$ is {\em finite} if the composite morphism 
\[
R \xrightarrow{\sigma} X \times_S X \xrightarrow{{\rm pr}_i} X
\]
with the $i$-th projection morphism ${\rm pr}_i$ is finite, for $i\in \{1,2\}$. 
\item 
Assume that $R$ and $X$ are reduced algebraic spaces over $S$. 
A relation $\sigma\colon R \to X\times_S X$ is 
a {\em set-theoretic equivalence relation} over  $S$
if, for every algebraically closed field $K$ and morphism $\Spec\,K \to S$, 
denoting $X_K:=X\times_S \Spec~K$ and $R_K:=R\times_S \Spec~K$, we have that
the image $\overline{R}_K(K)$ of the induced map 
\[
\sigma(K)\colon R_K(K) \to X_K(K) \times X_K(K)
\]
defines an equivalence relation on $X_K(K)$, i.e. 
the following hold: 
\begin{enumerate}
\item 
If $x \in X_K(K)$, then $(x, x) \in \overline{R}_K(K)$. 
\item 
If $(x, x') \in \overline{R}_K(K)$ with $x, x' \in X_K(K)$, then $(x', x) \in \overline{R}_K(K)$. 
\item 
If $(x, x'), (x', x'') \in \overline{R}_K(K)$ with $x, x', x'' \in X_K(K)$, then $(x, x'') \in \overline{R}_K(K)$. 
\end{enumerate}
\item Let $\sigma \colon R\to X\times_S X$ be a  set-theoretic equivalence relation. 
A {\em categorical quotient} of $X$ by $R$ over  $S$ is 
an $S$-morphism $q\colon X\to Y$ of algebraic spaces over $S$ such that 
\begin{enumerate}
\item $q\circ {\rm pr_1}\circ \sigma = q\circ {\rm pr_2}\circ \sigma$, and 
\item $Y$ is universal with respect to property (a), i.e. given any $S$-morphism $q'\colon X\to Y'$ to an algebraic space $Y'$ over $S$ such that $q'\circ {\rm pr_1}\circ \sigma = q'\circ {\rm pr_2}\circ \sigma$, 
there is a unique $S$-morphism $\pi\colon Y\to Y'$
satisfying $q'=\pi\circ q$. 
\end{enumerate}
\item Let $\sigma \colon R\to X\times_S X$ be a  set-theoretic equivalence relation. A categoric quotient $q\colon X\to Y$ of $X$ by $R$ is called a {\em geometric quotient} if 
$q$ is finite and the induced morphism $R \to (X\times_Y X)_{\red}$ is an isomorphism.
In this case, we denote $Y$ by $X/R$. 
\end{enumerate}
\end{dfn}
When no confusion arises, we will simply call a relation (resp. set-theoretic equivalence relation, \dots) on $X$ over $S$ as a relation (resp. set-theoretic equivalence relation, \dots) on $X$. 

\begin{ex}
Let $S$ be a scheme and let $f\colon X \to Y$ be an $S$-morphism of reduced algebraic spaces separated over $S$.

Then the induced closed immersion  
\[
(X \times_Y X)_{\red} \to X \times_S X
\]
defines a set-theoretic equivalence relation.\end{ex}

The following theorem is due to Koll\'ar: 

\begin{thm}\label{t-Kollar} 
Let $S$ be a noetherian $\mathbb F_p$-scheme and let $X$ be an algebraic space which is proper over $S$. Let $\sigma \colon R\to X\times_S X$ be a finite, set-theoretic equivalence relation. 

Then a geometric quotient $X\to X/R$ exists. 
\end{thm}

\begin{proof}
See \cite[Theorem 6]{kollar12}. 
\end{proof}

\begin{dfn}\label{d_l-equivalence}
Let $K$ be an algebraically closed field and let $X$ be an algebraic space which is proper over $\Spec\,K$. 
Let $L$ be a nef invertible sheaf on $X$.  
The $L$-{\em equivalence relation} on $X$ 
is the subset $R_L$ of $X(K) \times X(K)$ such that, for any $(x_1, x_2) \in X(K) \times X(K)$, we have that 
$(x_1, x_2) \in R_L$ if and only if 
there exists a morphism $j\colon C \to X$ 
from a one-dimensional proper connected $K$-scheme $C$ 
such that $x_1, x_2 \in j(C(K))$ and $j^*L$ is numerically trivial. 
Given a positive integer $m_0$, we say that the $L$-equivalence relation is {\em bounded} by $m_0$ if, in the notation above, 
we can always choose $C$ so that the number of irreducible components of $C$ is at most $m_0$. 
\end{dfn}

\begin{rem}
Note that, in general, the $L$-equivalence relation is not a set-theoretic equivalence relation. We refer to \cite[\S 5]{keel03} for an example of a nef invertible sheaf $L$ on a projective normal variety $X$ such that 
the $L$-equivalence relation is not bounded by any positive integer. 
\end{rem}

We now prove Theorem \ref{t-Keel-EWM}.

\begin{proof}[Proof of Theorem \ref{t-Keel-EWM}]
The only-if part is clear. 
We show the other implication. 
Assume that (1) and (2) hold. 
Let 
$g\colon Y \to Z$
be the morphism over $S$ induced by $L|_Y$. We may assume that $g_*\MO_Y=\MO_Z$. 
We have two set-theoretic equivalence relations on $Y \times_S Y$, given by 
\[
(Y \times_X Y)_{\red} \qquad\text{and}\qquad (Y \times_Z Y)_{\red}.
\]
We divide the proof into three steps. 
\setcounter{step}{0}
\begin{step}
In this step, 
we inductively define a reduced closed subspace $R^m$ of $Y \times_S Y$ 
for any  $m\in \Z_{\geq 1}$. 

We first set 
\[
R^1:=(Y \times_X Y) \cup  (Y \times_Z Y),
\]
equipped with the reduced closed subspace structure \cite[Definition II.6.9 and Proposition II.6.10]{knutson71}.
Assume that $R^m$ is already defined. 
Then we define $R^{m+1}$ as the image of the composite morphism 
\[
R^m_{12} \times_{Y_2} R^m_{23} \hookrightarrow 
(Y_1 \times_S Y_2) \times_{Y_2} (Y_2 \times_S Y_3) \to Y_1 \times_S Y_3
\]
equipped with the reduced subspace structure, 
where $Y_1, Y_2, Y_3:=Y$ and $R^m_{12}, R^m_{23}:=R^m$ are equipped with 
the corresponding projection morphisms. 
Since each $R^m$ contains the diagonal $\Delta_{Y/S}$ of $Y \times_S Y$, 
we have that 
\[
R^1 \subset R^2 \subset \cdots.
\]
\end{step}

\begin{step}\label{s-eq-termination}
Let $m_1=2m_0$. Then 
\[
R^{m}=R^{m_1}\qquad \text{for any  $m \geq m_1$.} 
\]
Thus, we define $R:=R^{m_1}$. 
\end{step}

\begin{proof}[Proof of Step~\ref{s-eq-termination}]
Let $K$ be an algebraically closed field. 
It is enough to show that $R^{m}=R^{m_1}$ for any $m \geq m_1$, under the assumption that $S=\Spec\, K$. 

Let $m\ge m_1$ and pick  $(y, \widetilde y)\in R^m(K)$. 
Let $x, \widetilde x \in X(K)$ be the images of $y$ and $\widetilde y$ respectively. 
Then there exist $\ell_0 \in \{1, 2, \cdots, m_0\}$ and 
$K$-curves $C_1, \cdots, C_{\ell_0}$ in $X$ 
such that $x \in C_1, \widetilde x \in C_{\ell_0}$ and $\bigcup_{i=1}^{\ell_0} C_i$ 
is connected. 
After possibly removing superfluous curves, 
we may assume that $C_i \cap C_{i+1}$ is not empty for each $i$. Let   $x_{i, i+1}\in X(K)$
be one of the intersection points.  
Let $C'_1, \dots, C'_{\ell_0}$ be $K$-curves in $Y$ such that 
$f(C'_i)=C_i$, $y \in C'_1$ and $\widetilde y \in C'_{\ell_0}$.
Let $y^{(i)}_{i, i+1}$ (resp.  $y^{(i+1)}_{i, i+1}$) be a closed point 
of $C'_i$ (resp.  $C'_{i+1}$) lying over $x_{i, i+1}$. 
Note that $L|_Y\cdot C'_i=0$ for all $i$. 
Since, for each $i$, we have 
\[
(y, y^{(1)}_{1, 2}) \in R^1, \qquad
(y^{(\ell_0)}_{\ell_0-1, \ell_0}, \widetilde y) \in R^1,
\]
\[
(y^{(i)}_{i-1, i}, y^{(i)}_{i, i+1}) \in R^1,
 \qquad \text{and}\qquad 
(y^{(i)}_{i, i+1}, y^{(i+1)}_{i, i+1}) \in R^1,
\]
it follows that $(y,\widetilde y)\in R^{m_1}$, as claimed. 
\end{proof}

\begin{step}\label{step-Keel-EWM}
We now prove Theorem~\ref{t-Keel-EWM}.

Let $R_Z$ be the image of $R$ in $Z \times_S Z$, 
equipped with the reduced closed subspace structure. 
By Step~\ref{s-eq-termination}, $R$ is a set-theoretic equivalence relation on $X$. 
Since $R$ contains $(Y \times_Z Y)_{\red}$, its image $R_Z$ is a set-theoretic equivalence relation on $Z$. 
Fix $i \in \{1,2\}$. 
We consider the commutative diagram: 
\[
\begin{CD}
R @>\widetilde{g}>> R_Z\\
@VVjV @VVj'V\\
Y \times_S Y @>g \times g >> Z \times_S Z\\
@VV{\rm pr}_iV @VV{\rm pr}'_iV\\
Y @>g>> Z,
\end{CD}
\]
where 
the upper vertical arrows are the induced closed immersions and 
the lower vertical arrows are the $i$-th projection morphisms. 

We now show that the induced morphism $\pi'_i:={\rm pr}'_i\circ j'\colon R_Z\to Z$ is finite, for $i\in \{1,2\}$. 
As $\pi'_i$ is proper, being finite is equivalent to being quasi-finite, 
i.e. it is enough to show that  all fibres are zero-dimensional. 
Therefore, we are reduced to consider
the case where $S=\Spec\,K$ for  an algebraically closed field $K$.  
We assume by contradiction  that $\pi'_i$ is not finite. Then 
there exists a closed point $z$ of $Z$ such that 
the fibre $R_{Z, z}$ of $\pi'_i$ over $z$ contains a $K$-curve $C$. 
Since $\widetilde g:R \to R_Z$ is a proper surjective morphism, 
there exist a closed point $y\in Y$ and a $K$-curve $C_Y$ in $R$ 
such that $\widetilde g(C_Y)=C$, ${\rm pr}_i \circ j(C_Y)=\{y\}$ 
and $g(y)=z$. 
The image $\overline{C}_Y$ of $C_Y$ by the other projection: 
${\rm pr}_{3-i} \circ j:R \to Y$ is a $K$-curve in $Y$ such that $L|_Y \cdot \overline{C}_Y=0$ and $g(\overline{C}_Y)$ is not a point. 
However, this contradicts the fact that $g$ is the morphism induced by $L|_Y$. 
Thus, $\pi'_i$ is finite as claimed. 
In particular, $R_Z$ is a finite set-theoretic equivalence relation on $Z$.

By Theorem\,\ref{t-Kollar}, there exists a geometric quotient $q\colon Z\to Z/R_Z$. 
In particular, $q$ is a morphism over $S$. 
Since
\[
q \circ {\rm pr}_1' \circ j'= q \circ {\rm pr}_2' \circ j',
\]
it follows from the diagram above that 
\begin{equation}\label{e-Keel-EWM1}
q \circ g \circ {\rm pr}_1 \circ j= q \circ g \circ {\rm pr}_2 \circ j.
\end{equation}

Since $f\colon Y\to X$ is finite, \cite[Example 5]{kollar12} implies that the geometric quotient $q'\colon Y\to W:=Y/(Y\times_X Y)_{\red}$ exists and there is a finite  universal homeomorphism $\sigma\colon W\to X$ such that $f=\sigma\circ q'$. In particular, $q'$ is a finite morphism. Since 
$R$ contains $(Y\times_X Y)_{\red}$, 
it follows from (\ref{e-Keel-EWM1}) and by (4) of Definition \ref{d-quotients} 
that the morphism $q \circ g\colon Y\to Z/R_Z$ uniquely factors through $W$. 
Let $h\colon W\to Z/R_Z$ be the induced $S$-morphism.

We now show that $L|_W$ is EWM over $S$. 
Let $s\in S$ be a  point and let  $V$ be an irreducible closed subspace of $W_s$. 
It is enough to show that $\dim h(V)<\dim V$ if and only if 
$(L|_W)^{\dim V}\cdot V=0$ 
(cf. Subsection \ref{s_properties-invertible-sheaves}). 
Let  $V'$ be an irreducible closed subspace of $Y_s$ such that $q'(V')=V$. Note that, since $q'$ is finite, we have that $\dim V'=\dim V$ and 
$(L|_W)^{\dim V}\cdot V=0$ if and only if 
$(L|_Y)^{\dim V'} \cdot V'=0$. 
Since $g$ is the morphism over $S$ induced by $L|_Y$, it follows that 
$(L|_Y)^{\dim V'} \cdot V'=0$ 
if and only if $\dim g(V')<\dim V'$. Since $q$ is a finite morphism, 
we have that $\dim h(V)=\dim q(g(V')) =\dim g(V')$, as claimed. Thus, $L|_W$ is EWM over $S$ and Lemma~\ref{l_EWM} implies that $L$ is EWM over $S$. 
 \qedhere
\end{step}
\end{proof}

\begin{thm}\label{t_AB-to-C}
Fix a positive integer $n$.

Then (Theorem~\ref{t-A})$_n$ and (Theorem~\ref{t-nume-triv4})$_n$ imply (Theorem~\ref{t-C})$_n$.
\end{thm}

\begin{proof}
Let $f:X \to S$ and $L$ be as in (Theorem~\ref{t-C})$_n$. 
By Lemma~\ref{l_reduced}, we may assume that $X$ and $S$ are reduced. 
 (Theorem~\ref{t-A})$_n$ implies that the restriction of $L$ to the normalisation $X^N$ of $X$
is semi-ample over $S$. 

\begin{claim}
There exists a positive integer $m_0$ such that, for all $s\in S$, the  $L|_{X_s}$-equivalence is bounded by $m_0$, 
\end{claim}

\begin{proof}[Proof of  Claim]
By assumption, for any point $s\in S$, we have that $L|_{X_s}$ is semi-ample. Let $g_s\colon X_s\to Z_s$ be the induced morphism. Let $m_s$ be the maximum number of irreducible components of any fibre of $g_s$. Then the $L|_{X_s}$-equivalence is bounded by $m_s$. 
Spreading $g_{\xi}$ out for any generic point $\xi$ of $S$, 
there exists an open dense subset $S^0$ of $S$ and morphisms 
\[
f^0\colon X^0:=f^{-1}(S^0) \overset{g^0}\to Z^0 \overset{h^0}\to S^0
\]
such that $f^0=f|_{f^{-1}(S^0)}$ and $g^0|_{X_s}=g_s$ for all $s\in S_0$. Thus, there exists a positive integer $m_1$ such that $m_s\le m_1$ for all $s\in S^0$. By noetherian induction, we may find a positive integer $m_2$ such that $m_s\le m_2$ for all $s\in S\setminus S^0$. Thus, it is enough to take $m_0=\max\{m_1,m_2\}$.
\end{proof}
Thus, $L$ is EWM over $S$ by Theorem~\ref{t-Keel-EWM}. 
Let 
\[
f\colon X \xrightarrow{g} Z \to S,
\]
be the morphisms induced by $L$. 
Lemma \ref{l-exc-fin-sub} implies that $Z$ is excellent. 

By Theorem~\ref{t-nume-triv-as},
there exists an invertible sheaf $L_Z$ on $Z$ and a positive integer $m$ such that $L^{\otimes m}=g^*L_Z$. 
Since $g$ contracts the $L$-trivial curves, $L_Z$ is ample over $S$ 
by the Nakai--Moishezon criterion (cf. \cite[Theorem 3.11]{kollar90}, \cite[Proposition 1.41]{KM98}).  
In particular,  $L$ is semi-ample over $S$. 
\end{proof}

 \section{Proof of the main theorems}

\begin{proof}[Proof of Theorem~\ref{t_main}]
By Theorem \ref{t_C-to-A}, Theorem \ref{t_A-to-B} and Theorem \ref{t_AB-to-C}, 
(Theorem \ref{t-C}$)_n$ holds for any $n \in \Z_{\geq 0}$. 
Therefore Theorem~\ref{t_main} holds if $X$ is finite dimensional. 
By Remark \ref{r-reduce-hensel}, 
we can reduce the general case 
to this case, after possibly replacing $S$ by the affine spectrum of a stalk. 
\end{proof}

\begin{lemma}\label{l_uncountable-field}
Let $k$ be an uncountable field and 
let $f\colon X\to S$ be a projective $k$-morphism of schemes of finite type over $k$.  
Let $L$ be an invertible sheaf on $X$ such that  $L|_{X_s}$ is semi-ample for all the closed points $s\in S$. 

Then $L|_{X_s}$ is semi-ample 
for any point $s\in S$.
\end{lemma}

\begin{proof} 
We show the lemma by induction on the dimension of $S$. 
If $\dim S=0$, then the claim is clear. 
Thus, we may assume that $\dim S>0$ and that the claim holds if the dimension 
of the base is smaller than $\dim S$. 
In particular, it is enough to show that $L|_{X_\xi}$ is semi-ample 
for the generic point $\xi \in S$ of an irreducible component of $S$. 
Replacing $S$ by an open neighbourhood of $\xi$, 
we are reduced to the case where 
$S$ is an affine integral scheme such that $f$ is flat. 

By
the semicontinuity theorem 
\cite[Theorem III.12.8]{hartshorne77}, 
for any positive integer $m$, 
 there exist a positive integer $c_m$ and 
a non-empty affine open subset $U_m \subset S$ such that 
\[
c_m=\dim_{k(s)} H^0(X_s,L^{\otimes m}|_{X_s})
\]
for any point $s \in U_m$. Since $k$ is uncountable, 
there exists a closed point 
\[
t\in \bigcap_{m\in \Z_{>0}} U_m.
\]
As $L|_{X_t}$ is semi-ample, 
there exists a positive integer $m_0$ such that 
$L^{\otimes m_0}|_{X_t}$ is globally generated. 
By Grauert's theorem \cite[Corollary III.12.9]{hartshorne77}, 
the restriction map 
\[
H^0(f^{-1}(U_{m_0}), L^{\otimes m_0}|_{f^{-1}(U_{m_0})}) \to H^0(X_t, L^{\otimes m_0}|_{X_t})
\]
is surjective. 
Since the base locus of the linear system associated to $L^{\otimes m_0}$ is a closed subset of $X$, 
it is disjoint from $X_t$. 
In particular, $L|_{X_{\xi}}$ is semi-ample, as desired. 
\end{proof}

\begin{proof}[Proof of Theorem~\ref{t_main2}]
Theorem~\ref{t_main} and Lemma~\ref{l_uncountable-field} immediately imply the claim. 
\end{proof}

\section{Examples}

\subsection{Examples over $\overline{\mathbb F}_p$}
The following example shows that, over countable fields, we need to consider not only closed points of $S$ but all the scheme-theoretic points of $S$ in Theorem~\ref{t_main} (cf. Theorem \ref{t_main2}). 

\begin{ex}\label{e_scheme}
Let $E$ be an elliptic curve over $\overline{\F}_p$. 
Let $X:=E \times E$ and $S:=E$. 
Let $f\colon X \to S$ be the first projection. 
Let $L:=\MO_X(\Delta-Z)$, where $\Delta$ is the diagonal divisor of $X=E \times E$ and 
$Z:=E \times \{Q\}$ for a closed point $Q \in E$. 
By \cite[Example 1.46]{KM98}, $L$ is $f$-nef but not $f$-semi-ample. 
Note that $L|_{X_s}$ is semi-ample for all the closed points $s \in S$ since the base field is $\overline{\F}_p$. 
On the other hand, Theorem~\ref{t_main} implies that  $L|_{X_{\xi}}$ is not semi-ample for the generic point $\xi $ of $S$.
\end{ex}

\subsection{Counterexamples in characteristic zero}\label{s_c-zero}

The goal of this subsection is to show that Theorem~\ref{t_main} does not hold in characteristic zero. 
The following result is due to Keel: 

\begin{prop}\label{p_keel} 
Let $k$ be an algebraically closed field of characteristic zero. 
Let $C$ be a smooth projective curve over $k$ and whose genus is at least two. 
Let $X:=C \times_k C$ and let $\pi_i\colon X\to C$ be the $i$-th projection 
for  $i \in \{1,2\}$. Let 
 $\Delta \subset X$ be the diagonal and let 
\[
L:=\MO_X(K_X-\pi_1^*K_C+\Delta).
\]

Then the following hold:
\begin{enumerate}
\item $L$ is nef and big. 
\item $L|_{\Delta} \simeq \MO_{\Delta}$ and $L \cdot D>0$ for a curve $D$ in $X$ other than $\Delta$. 
\item $L|_{2\Delta}$ is not semi-ample. 
\end{enumerate}
\end{prop}

\begin{proof}
By  \cite[Theorem 3.0]{keel99}, (1) holds.  \cite[Lemma 3.2]{keel99} implies  (2) and  \cite[Lemma 3.4]{keel99} implies (3).
\end{proof}

\begin{ex}\label{e_scheme2}
Let $k$ be an algebraically closed field of characteristic zero. 
Let $C$ be a smooth projective curve over $k$ 
such that the genus of $C$ is at least three and $C$ is not hyperelliptic. 
Let $X=C\times_k C$ and let $\Delta\subset X$ be the diagonal. Let $L$ be as in  Proposition~\ref{p_keel}. 
Then, by \cite[Exercise V.D-2]{ACGH85}, there exists a birational morphism $f\colon X\to S$ onto a projective surface $S$ such that the exceptional locus of $f$ is $\Delta$. Moreover, if $s_0=f(\Delta)$, then \cite[Exercise VI.A-5]{ACGH85} implies that $X_{s_0}=\Delta$, i.e. $X_{s_0}$ is reduced. 
Thus, (2)  of Proposition~\ref{p_keel} implies that   $L|_{X_s}$ is semi-ample for all $s\in S$. However (3)  of Proposition~\ref{p_keel} implies that $L$ is not $f$-semi-ample. 
\end{ex}

\section{Errata}

This section is not contained in the published version. 
The proof of Theorem \ref{t-alteration} contains a gap, which was kindly
pointed out to us by Adrian Langer. The issue affects the proof of
Proposition \ref{p-normal-nt1}, where Theorem \ref{t-alteration} is used.
We now give a new proof of Proposition \ref{p-normal-nt1} which does not
rely on Theorem \ref{t-alteration}.

\begin{prop}[=Proposition \ref{p-normal-nt1}]\label{p errata}
Fix a positive integer $n$ and assume 
$({\rm Theorem}~{\rm C})_{n-1}$. 
Let $f\colon X \to S$ be a proper morphism 
of excellent $\F_p$-schemes satisfying $f_*\MO_X=\MO_S$, 
where $X$ is a normal scheme of dimension $n$. 
Let $L$ be an invertible sheaf on $X$ such that 
$L|_{X_s}$ is semi-ample for all the points $s \in S$ and 
$L|_{X_\xi}$ is numerically trivial for all the generic points $\xi$ of $S$.

Then $L$ is $f$-semi-ample. 
\end{prop}

\setcounter{step}{0}
\begin{proof}
The proof consists of four steps.

\begin{step}\label{s1 errata}
In order to prove Proposition \ref{p errata}, 
we may assume that 
\begin{enumerate}
\item $f \colon X \to S$ is projective, 
\item the generic fibre $X_\xi$ is geometrically integral and geometrically normal, and 
\item $S = \Spec\,R$ for a complete noetherian normal local ring $(R, \m)$. 
\end{enumerate}
\end{step}

\begin{proof}[Proof of Step \ref{s1 errata}]
By Chow's lemma, we may assume (1) (Lemma \ref{l_pullback}(4)). 
For the algebraic closure $\overline{\xi}$ of $\xi$, 
the geometric generic fibre $X_{\overline{\xi}}$ of $f$ is an irreducible projective scheme over $\overline{\xi}$ 
(cf.\ \cite[Lemma 2.2]{Tan18}). Then we can find a commutative diagram 
\[
\begin{tikzcd}
X' \arrow[r, "\alpha"] \arrow[d, "f'"] & 
X \arrow[d, "f"]\\
S' = \Spec\,R' \arrow[r, "\beta"]  & 
S=\Spec\,R.
\end{tikzcd}
\]
of projective surjective morphisms such that 
\begin{itemize}
\item $X'$ is a normal scheme of dimension $n$, 
\item $R'$ is an integral normal domain, 
\item $f'_*\MO_{X'} = \MO_{S'}$, 
\item each of $\alpha$ and $\beta$ is a finite morphism, and 
\item the geometric generic fibre of $f'$ is a normal integral scheme.
\end{itemize}
Replacing $f' \colon  X' \to S'$ by  $f \colon  X \to S$, 
we may assume (2) (Lemma \ref{l-stein}(2), Lemma \ref{l_pullback}(4)). 
Taking a localisation, we may assume that $S = \Spec\,R$ for a normal local ring $(R, \m)$ (Remark \ref{r-reduce-hensel}). 
For the completion $\widehat{R}$ of $R$, 
we see that $X \times_{\Spec\,R} \Spec\,\widehat{R}$ is normal. 
Again by Remark \ref{r-reduce-hensel}, 
we may assume (3). 
This completes the proof of Step \ref{s1 errata}. 
\qedhere

\end{proof}

\begin{step}\label{s2 errata}
Set $k:=R/\m$. 
Then there exist 
an intermediate ring $k \subset R_1 \subset R$, 
a cartesian diagram 
\[
\begin{tikzcd}
X \arrow[r, "\alpha_1"] \arrow[d, "f"] & 
X_1 \arrow[d, "f_1"]\\
S= \Spec\,R \arrow[r, "\beta_1"]  & 
S_1 := \Spec\,R_1
\end{tikzcd}
\]
consisting of projective morphisms, 
and an invertible sheaf $L_1$ on $X_1$ such that 
\begin{enumerate}
\item $R_1$ is a finitely generated $k$-algebra, 
\item $X_1$ is an integral scheme, 
\item $\m_1 := \m \cap R_1$ is a maximal ideal satisfying $k \xrightarrow{\simeq} R_1/\m_1$, and 
\item $\alpha_1^*L_1 \simeq L$. 
\end{enumerate}
\end{step}

\begin{proof}[Proof of Step \ref{s2 errata}]
As $(R, \m)$ is a complete noetherian local ring, 
$R$ is a $k$-algebra with $k \xrightarrow{\simeq} R/\m$. 
We have $X = \Proj\,R[y_0, ..., y_N]/J$ for some homogeneous ideal $J$ 
of $R[y_0, ..., y_N]$. 
Hence there exists 
a finitely generated $k$-subalgebra $R_1$ of $R$ and a homogeneous ideal $J_1$ 
of $R[y_0, ..., y_N]$ satisfying $X_1 \times_{\Spec R_1} \Spec R = X$ 
for $X_1 := \Proj R_1[y_0, ..., y_N]/J_1$. 
Thus (1) holds. 
Enlarging $R_1$ if necessary, 
we can find an invertible sheaf $L_1$ on $X_1$ 
satisfying (4). Since $\m_1 := \m \cap R_1$ is a prime ideal of $R_1$ satisfying 
${\rm id} : k \to R_1/\m_1 \hookrightarrow R/\m =k$, 
we get $k \xrightarrow{\simeq} R_1/\m_1$. 
Thus (3) holds.

It is enough to reduce the problem to the case when (2) holds. 
We now prove that we may assume that $X_1$ is reduced. 
We have the induced morphisms 
\[
\theta: X = X_{\red} \xrightarrow{i_1} (X_1)_{\red} \times_{S_1} S 
\overset{i_2}{\hookrightarrow} X_1 \times_{S_1} S. 
\]
where the composition $\theta$ is an isomorphism and $i_2$ is a surjective closed immersion. 
In particular, $i_1$ is an affine morphism which is a homeomorphism. 
Then the corresponding composite ring homomorphism 
\[
\theta^{\sharp} : 
\MO_{X_1 \times_{S_1} S} \xrightarrow{i_2^{\sharp}}  \MO_{(X_1)_{\red} \times_{S_1} S}  \xrightarrow{i_1^{\sharp}}   \MO_X
\]
is an isomorphism, where $i_2^{\sharp}$ is surjective. 
Therefore, both $i_1^{\sharp}$ and $i_2^{\sharp}$ are isomorphisms, and hence so are $i_1$ and $i_2$. 
In what follows, we assume that $X_1$ is reduced.

Finally, it suffices to reduce the problem to the case when $X_1$ is an integral scheme. 
If we have $X_1 = Y_1 \cup Z_1$ for some 
reduced closed subschemes $Y_1$ and $Z_1$ of $X_1$, 
the irreducibility of $X$ enables us to deduce 
that $(Y_1 \times_{S_1} S)_{\red} = X$ or $(Z_1 \times_{S_1} S)_{\red} = X$. 
Repeating this procedure, 
we can find an irreducible component $Y_1$ of $X_1$, equipped with the reduced scheme structure, such that $(Y_1 \times_{S_1} S)_{\red} = X$. 
By the following inclusions of closed subschemes  of $X_1 \times_{S_1} S$:   
\[
X = (Y_1 \times_{S_1} S)_{\red} \subset
Y_1 \times_{S_1} S \subset X_1 \times_{S_1} S = X, 
\]
we get the scheme-theoretic equality $X = (Y_1 \times_{S_1} S)_{\red}$. 
Thus we may assume (2). 
This completes the proof of Step \ref{s2 errata}. 
\end{proof}

\begin{step}\label{s3 errata}
There exists a commutative diagram  
\begin{equation}\label{e1 s3 errata}
\begin{tikzcd}
Y_1 \arrow[r, "\alpha_1"] \arrow[d, "g_1"] & 
X_1 \arrow[d, "f_1"]\\
T_1 \arrow[r, "\beta_1"]  & 
S_1.
\end{tikzcd}
\end{equation}
consisting of projective morphisms of quasi-projective varieties over $k$ such that 
\begin{enumerate}
\item $Y_1$ is normal, $T_1$ is regular, 
\item $\alpha_1$ and $\beta_1$ are generically finite 
projective surjective morphisms, 
\item $g_1$ is equi-dimensional (i.e., 
for every $t \in T_1$, $g_1^{-1}(t)$ is equi-dimensional 
and $\dim g_1^{-1}(t)=\dim Y_1 - \dim T_1$), 
\item $(g_1)_*\MO_{Y_1} = \MO_{T_1}$, and 
\item 
the diagram 
\[
\begin{tikzcd}
Y_1 \times_{S_1} \xi_1 \arrow[r, "\alpha_1"] \arrow[d, "g_1"] & 
X_1 \times_{S_1} \xi_1 \arrow[d, "f_1"]\\
T_1 \times_{S_1} \xi_1 \arrow[r, "\beta_1"]  & 
\xi_1
\end{tikzcd}
\]
obtained by applying the base change $(-) \times_{S_1} \xi_1$ to 
the diagram (\ref{e1 s3 errata}) is cartesian, 
where $\xi_1$ denotes the generic point of $S_1$. 
\end{enumerate}
\end{step}

\begin{proof}[Proof of Step \ref{s3 errata}]
Recall that $S_1=\Spec\,R_1$ is an affine variety over $k$ and 
$X_1$ is a quasi-projective variety over $k$ (Step \ref{s2 errata}). 
We now construct a commutative diagram \[
\begin{tikzcd}
Y_1 := X'''_1 \arrow[r] \arrow[d, "f'''_1"] & 
X''_1 \arrow[r] \arrow[d, "f''_1"] & 
X'_1 \arrow[r] \arrow[d, "f'_1"] & 
X_1 \arrow[d, "f_1"]\\
T_1 := S'''_1 \arrow[r] & 
S''_1 \arrow[r] & 
S'_1 \arrow[r]  & 
S_1.
\end{tikzcd}
\]
By \cite[Th\'eor\`eme 5.2.2]{RG71}, there exists a 
flattening $f'_1 : X'_1 \to S'_1$ of $f_1: X_1 \to S_1$, so that 
\begin{itemize}
\item 
$X'_1$ and $S'_1$ are quasi-projective varieties over $k$, 
\item 
$X'_1 \to X_1$ and $S'_1 \to S_1$ are projective birational morphisms, 
\item  $f'_1 : X'_1 \to S'_1$ is a flat projective morphism, and 
\item $(X'_1)_{\xi'_1} \xrightarrow{\simeq} (X_1)_{\xi_1}$, 
where $\xi'_1$ denotes the generic point of $S'_1$. 
\end{itemize}
Let $X''_1$ be the normalisation of $X'_1$ and 
let $S''_1$ be the Stein factorisation of the induced composition 
$X''_1 \to X'_1 \to S'_1$. 
In particular, $(f''_1)_*\MO_{X''_1}= \MO_{S''_1}$. 
Take an alteration $S'''_1 \to S''_1$ of $S''_1$ \cite[Theorem 4.1]{jong96} (i.e., 
a generically finite projective morphism from a regular integral scheme $S'''_1$). 
Since $(X_1)_{\xi_1} \times_{\xi_1} \xi \simeq X_{\xi}$ 
(Step \ref{s2 errata}) and $X_{\xi}$ is geometrically integral and geometrically normal (Step \ref{s1 errata}(2)),  
so is the generic fibre $(X_1)_{\xi_1}$ of $f_1$. 
By construction, we have $(X''_1)_{\xi''_1} \xrightarrow{\simeq} (X'_1)_{\xi'_1} \xrightarrow{\simeq}(X_1)_{\xi_1}$ 
for the generic point $\xi''_1$ of $S''_1$. 
Then the generic fibre $(X''_1)_{\xi''_1}$ 
of $f''_1 : X''_1 \to S''_1$ is 
geometrically integral and geometrically normal, and hence so is  the generic fibre of the base change $X''_1 \times_{S''_1} S'''_1 \to S'''_1$.  
In particular, there exists a unique irreducible component $W$ of $X''_1 \times_{S''_1} S'''_1$ dominating $S'''_1$. 
Let $X'''_1$ be the normalisation of $W$. 
Then (1)--(3) hold. 
Note that 
 $X''_1 \times_{S''_1} S'''_1 \to S'''_1$ and 
$X'''_1 \to S'''_1$ coincide over the generic point $\xi'''_1 \in S'''_1$. 
Thus (5) holds. 
Moreover, the ring homomorphism $\MO_{S'''_1} \to (f'''_1)_*\MO_{X'''_1}$ is an isomorphism at the generic point. 
This, together with the normality of $S'''_1$,  implies 
that $\MO_{S'''_1} \xrightarrow{\simeq} (f'''_1)_*\MO_{X'''_1}$. 
Thus (4) holds. 
This completes the proof of Step \ref{s3 errata}. 
\end{proof}

\begin{step}\label{s4 errata}
$L$ is $f$-semi-ample. 
\end{step}

\begin{proof}[Proof of Step \ref{s4 errata}]
Set $R_2 := (R_1)_{\m_1}$, $\m_2 := \m_1 (R_1)_{\m_1}$, and 
$S_2 := \Spec\,R_2$. 
By $R_1 \subset R_2 \subset R$, we have $S \to S_2 \to S_1$. 
Applying  the base changes $(-) \times_{S_1} S_2$ and $(-) \times_{S_1} S$ to (\ref{e1 s3 errata}), we obtain commutative diagrams 
\[
\begin{tikzcd}
Y_2 \arrow[r, "\alpha_2"] \arrow[d, "g_2"] & 
X_2 \arrow[d, "f_2"]\\
T_2 \arrow[r, "\beta_2"]  & 
S_2.
\end{tikzcd}
\qquad\qquad 
\begin{tikzcd}
Y \arrow[r, "\alpha"] \arrow[d, "g"] & 
X \arrow[d, "f"]\\
T \arrow[r, "\beta"]  & 
S.
\end{tikzcd}
\]
By Lemma \ref{l-nef-special}, the pullback $L_{X_2}$ of $L_{X_1}$ on $X_2$ is nef over $S_2$, and hence so is its pullback $L_{Y_2}$ on $Y_2$. 
As $(-) \times_{S_1} S_2$ is a localisation, 
$T_2$ is regular and we have $\MO_{T_2} \xrightarrow{\simeq} (g_2)_*\MO_{Y_2}$. 
After replacing $L_{X_1}$ by a multiple if necessary, 
we get $L_{Y_2} \simeq g_2^*L_{T_2}$ for some $L_{T_2} \in \Pic\,T_2$ (Lemma \ref{l-flat-descent}). 

Then it is enough to show that 
the pullback $L_{T}$ of $L_{T_2}$ is $\beta$-semi-ample (Lemma \ref{l_pullback}(4)). 
Since $\beta_1 : T_1 \to S_1$ is a generically finite projective morphism (Step \ref{s3 errata}(2)), 
$\beta$ is a projective morphism such that 
$T_\xi \to \xi$ is a finite morphism. 
In particular, $L_{T}$ is $\beta$-big and 
we have $\mathbb E_{\beta}(L_T) \subset \beta^{-1}(\overline{S})$ 
for some proper closed subset $\overline{S}$ of $S$. 
We equip $\overline{S}$ with the reduced scheme structure. 
By Proposition \ref{p-relative-keel}, 
it suffices to show that $L_{T}|_{\beta^{-1}(\overline{S})}$ is semi-ample over $S$. 
Since $g_1$ has connected fibres (Step \ref{s3 errata}(4)), 
so does a base change $g$. 
Then it is enough to prove that the pullback 
$L_Y|_{g^{-1}(\beta^{-1}(\overline{S}))}$
of 
$L_{T}|_{\beta^{-1}(\overline{S})}$ 
is semi-ample over $S$, which follows from $({\rm Theorem}~{\rm C})_{n-1}$. 
This completes the proof of Step \ref{s4 errata}. 
\end{proof}
Step \ref{s4 errata} completes the proof of 
Proposition \ref{p errata}. 
\end{proof}

\bibliographystyle{amsalpha}
\bibliography{Library}

\end{document}